\documentclass[reqno]{amsart}
\usepackage{dashint}
\usepackage{amssymb,amsmath,amsbsy,array,latexsym,hyperref}

\theoremstyle{plain}
\newtheorem{theorem}                 {Theorem}   [section]
\newtheorem{corollary}   [theorem]   {Corollary}
\newtheorem{lemma}       [theorem]   {Lemma}
\newtheorem{proposition} [theorem]   {Proposition}
\newtheorem{definition}  [theorem]   {Definition}

\theoremstyle{remark}

\newtheorem{remark}      [theorem]   {Remark}
\newtheorem*{notations}              {\textbf{Notation}}
\numberwithin{equation} {section}
\begin{document}
\title[Sobolev classes and horizontal energy minimizers]
{Sobolev classes and horizontal energy minimizers between
Carnot-Carath\'{e}odory spaces}

\author{Kanghai Tan}

\address{Department of Applied Mathematics\\
Nanjing University of Science and Technology\\ 210094, Nanjing,
The People's Republic of China}
\email{tankanghai2000@yahoo.com.cn}
\begin{abstract}
The notion of horizontal energy minimizers between C-C spaces is
introduced. We prove existence of such energy minimizers when the
domain is a $C^{2}$, noncharacteristic bounded open set in a C-C
space and the target is a C-C space of Carnot type. \vskip.2truecm
{\small\noindent {\bf Keywords:} Sobolev mappings,
Carnot-Carath\'{e}odory spaces, Carnot groups, energy minimizers,
existence, regularity
\vskip.2truecm
\noindent {\bf
2000Mathematics Subject Classification:} 58E20, 49N60, 53C22.}
\end{abstract}
\maketitle

\section{Introduction}\label{intro}
Recently many people paid their attentions to the study of
analysis and geometry in metric measure spaces in particular in
Carnot-Carath\'{e}odory (written as C-C for brevity) spaces, see
\cite{HK}, \cite{Sub-Rie} and references therein. In this
direction, there has been a number of works devoted to the notions
of Sobolev functions and mappings on metric spaces. Let us in
particular mention the definitions proposed by Korevaar-Schoen in
\cite{KS}, by Haj\l asz in \cite{Ha} and by Reshetnyak in
\cite{Re}, see \cite{Ch}, \cite{Sh}, \cite{FLW1}, \cite{LLW} and
\cite{LP} for other definitions and generalizations. The Sobolev
spaces of \cite{Ch}, \cite{Sh} and \cite{Ha} are originally
defined in metric measure spaces for real-valued functions. These
classes of Sobolev functions are equivalent as sets when the
Sobolev exponent is larger than one, and all equivalent to the
horizontal Sobolev spaces (\cite{GN}) when the domain is a C-C
space satisfying suitable conditions, see \cite{FHK}, \cite{Sh}
and \cite{HK}. The notions in \cite{Ch}, \cite{Ha}, \cite{KS},
\cite{Re}, \cite{Sh} can be extended to define Sobolev mappings
between metric measure spaces in particular C-C spaces, see
\cite{HKST} and \cite{Vo1}.

On the other hand, in \cite{KS} Korevaar-Schoen used their
developed theory of Sobolev mappings to study harmonic mappings
from smooth Riemannian manifolds to nonpositive curvature spaces.
Let us briefly recall their ideas. Assume that $\Omega$ is a
smooth domain in $R^{n}$ and $M$ is a separable metric space with
a metric $d$. A function $u\in{L^{\alpha}}(\Omega, M)$ is in
$KS^{1,\alpha}(\Omega,M)$ if
$$
E^{\alpha}(u,\Omega)=\sup_{f\in{C_{c}(\Omega,[0,1])}}
\limsup_{\epsilon\rightarrow
0}\int_{\Omega}f(x)\dashint_{B_{\epsilon}(x)}
\left(\frac{d(u(x),u(y))}{\epsilon}\right)^{\alpha}dydx
$$
is finite where $C_{c}(\Omega,[0,1])$ is the set of all compactly
supported functions in $\Omega$ taking values in the interval
$[0,1]$. When $\Omega$ is an open set in a smooth Riemannian
manifold, the definition is similar. If
$u\in{KS^{1,\alpha}(\Omega,M)}$, $E^{\alpha}(u,\Omega)$ is called
the energy of the mapping $u$. Roughly speaking, the story of
\cite{KS} is based on a subpartitional lemma (\cite{KS}, Lemma
1.3.1). It follows from the subpartitional lemma that
$E^{\alpha}(u,\Omega)$ is lower semicontinuous with respect to the
topology of $L^{\alpha}(\Omega, M)$ and $KS^{1,\alpha}(\Omega,M)$
possesses some type of precompactness property(\cite{KS}, Theorem
1.13). Korevaar-Schoen proved a satisfactory existence and
regularity theory for energy minimizers of $E^{\alpha}(u,\Omega)$
when the target is a nonpositive curvature space (in the sense of
Alexandrov). In \cite{EF} and \cite{Fu}, Eells and Fuglede made a
systematic generalization of the Korevaar-Schoen's results to
Riemannian polyhedra. For similar results but with different
methods we refer to \cite{Jo2}, \cite{Jo3}, \cite{Jo4}, and
\cite{Jo5}.

We briefly recall the definition of C-C spaces (or sub-Riemannian
manifolds), in particular of Carnot groups. Let $\Delta$ be a
smooth distribution in $R^{n}$ satisfying the H\"{o}rmander
condition and endowed with an inner product $<\cdot,\cdot>_{c}$.
The structure of $(\Delta,<\cdot,\cdot>_{c})$ yields the C-C
metric $d_{c}$, see Section \ref{sec:1-1} for details.
$(R^{n},\Delta,d_{c})$ is called a C-C space (if $R^{n}$ is
replaced by a smooth manifold $M$, $(M,\Delta,d_{c})$ is called a
sub-Riemannian manifold). Carnot groups are most interesting C-C
spaces. A \textit{Carnot group} $G$ is a connected, simply
connected Lie group whose Lie algebra $\mathcal{G}$ admits the
grading $\mathcal{G}=V_{1}\bigoplus\cdots\bigoplus V_{l}$, with
$[V_{1},V_{i}]=V_{i+1}$, for any $1\leq i\leq l-1$ and
$[V_{1},V_{l}]={0}$ (the integer $l$ is called the step of $G$).
Let $\{e_{1},\cdots,e_{n}\}$ be a basis of $\mathcal{G}$ with
$n=\sum_{i=1}^{l}\dim(V_{i})$. Let $X_{i}(g)=(L_{g})_{*}e_{i}$ for
$i=1,\cdots,k:=\dim(V_{1})$ where $(L_{g})_{*}$ is the
differential of the left translation
$L_{g}(g^{\prime})=gg^{\prime}$ and let
$Y_{i}(g)=(L_{g})_{*}e_{i+k}$ for $i=1,\cdots,n-k$. We call the
system of left-invariant vector fields
$\Delta:=V_{1}=\textrm{span}\{X_{1},\cdots,X_{k}\}$ the
\textit{horizontal bundle} of $G$. If we equip $\Delta$ an inner
product $<\cdot,\cdot>_{c}$ such that $\{X_{1},\cdots,X_{k}\}$ is
an orthonormal basis of $\Delta$, $(G,\Delta,<\cdot,\cdot>_{c})$
is an equiregular sub-Riemannian manifold. In
$(G,\Delta,<\cdot,\cdot>_{c})$, $d_{c}$ is invariant with respect
to left translation, that is $d_{c}(p_{0}p,p_{0}q)=d_{c}(p,q)$ for
any $p_{0},p,q\in{G}$, and is 1-homogeneous with respect to the
natural dilations, that is
$d_{c}(\delta_{s}p,\delta_{s}q)=sd_{c}(p,q)$ for any
$s>0,p,q\in{G}$, where
$\delta_{s}p=\exp(\sum_{i=1}^{l}s^{i}\xi_{i})$ for
$p=\exp(\sum_{i=1}^{l}\xi_{i}), \xi_{i}\in{V_{i}}$. We usually
identify $G$ with $R^{n}$ by the exponential map  and use
$(R^{n},V_{1},\delta_{\lambda})$ to denote $G$.
$Q=\sum_{i=1}^{l}i\dim(V_{i})$ is called homogeneous dimension of
$G$ and $\mathcal{L}^{n}$ is the Haar measure of $G$. It is easy
to prove that
\begin{equation}\label{sec1:coefficients}
X_{j}(x)=\frac{\partial}{\partial
x_{j}}+\sum_{i=k+1}^{n}a_{i}^{j}(x)\frac{\partial}{\partial
x_{j}},\quad X_{j}(0)=e_{j}=\frac{\partial}{\partial x_{j}},\quad
j=1,\cdots,k,
\end{equation}
where $a_{i}^{j}(x)=a_{i}^{j}(x_{1},\cdots,x_{k})$ are polynomials
such that
$a_{i}^{j}(\delta_{\lambda}x)=\lambda^{\alpha_{i}-\alpha_{j}}a_{i}^{j}(x)$.
The simplest noncommutative Carnot group is the Heisenberg group
$H^{m}$ which is, by definition, $R^{2m+1}$ with the group law
$pp^{\prime}=(z+z^{\prime},t+t^{\prime}+2\omega(z,z^{\prime}))$
where $p=(z,t),p^{\prime}=(z^{\prime},t^{\prime})\in{R^{2m}\times
R}$ and  $\omega$ stands for the standard symplectic form in
$R^{2n}$. For more about Carnot groups, see \cite{FS} and
\cite{Stein}.

 In this paper we want to generalize the theory
of harmonic mappings to C-C spaces in particular to Carnot groups.
In \cite{CL} Capogna and Lin made the first step in this
direction. Using the energy of Sobolev mappings of
Korevaar-Schoen, they considered energy minimizers with   smooth
Euclidean domain and target Heisenberg group  $H^{m}$ endowed with
a C-C metric. Note that Heisenberg group does not possess any
curvature bound in the sense of Alexandrov and the arguments in
\cite{KS} is not valid in this case. Capogna and Lin made full use
of the differential structure of the domain and the target to
characterize the Sobolev mappings and explicitly described the
energy. It turns out that these Sobolev mappings are weakly
contact (satisfying a Legendrian condition) while the energy is
not a Dirichlet integral (except the case when $\alpha =2$).
Precisely they proved the following:

\begin{theorem}\label{intro:thmcl}
Let $\Omega$ be a bounded domain in $R^{n}$ and $\alpha \geq 2$.
Then $u=(z,t)=(x,y,t)\in{KS^{1,\alpha}(\Omega,H^{m})}$ if and only
if $z\in{W^{1,\alpha}(\Omega)}$ and $t\in{L^{\frac{\alpha}{2}}}$
is weakly differentiable, and for a.e. $p\in{\Omega}$, $i=1,
\cdots,n$, $\partial_{p_{i}}
t=2(y\partial_{p_{i}}x-x\partial_{p_{i}}y)\in{L^{\beta}(\Omega)}$
with $\beta=\frac{n\alpha}{2n-\alpha}$. Moreover the energy can be
written as
\begin{equation}\label{intro:equcl}
E^{\alpha}(u,\Omega)=C\int_{\Omega}\int_{B_{1}}|\nabla z(p)\cdot
\omega|^{\alpha}d\omega dp.
\end{equation}
\end{theorem}
In general, the energy of Korevaar-Schoen can not be written as a
Dirichlet integral
\begin{equation}\label{intro:diri}
C\int_{\Omega}|\nabla \bar{u}(p)|^{\alpha}dp
\end{equation}
for some $\bar{u}$ related to $u$ when the target is not the real
line, see also \cite{KS}. It is easily seen that only when
$\alpha=2$, (\ref{intro:equcl}) has the form of
(\ref{intro:diri}). In our opinion in general the energy has the
form of (\ref{intro:diri}) is a necessary condition to make the
energy minimizing problem solvable when the target does not
possess any curvature bound. We remark that the method in
\cite{CL} used to characterize Sobolev mappings is not valid for
$1\leq\alpha<2$ due to the non-isotropic property of the gauge
distance.

To generalize the concept of harmonic mapping to C-C spaces we
must introduce a natural energy which not only has ``good" form
(like Dirichlet integral) but also inherits some essential nature
from the considered C-C spaces. To this end we first study the
energy of Korevaar-Schoen. We will show that the energy of
Korevaar-Schoen is not the one we expected. We will give an
explicit description of the energy of Korevaar-Schoen when both
the domain and the target are Carnot groups, see Theorem
\ref{thm:energychara}. That is,
\begin{equation}\label{intro:explicitcg}
E^{\alpha}(u,\Omega)=C\int_{\Omega}\int_{B_{c}(0,1)}\widetilde{\rho}(Du(p)(\omega))^{\alpha}
d\omega dp
\end{equation}
where $\Omega\subset G$ is a bounded open set of Carnot group $G$
with a homogeneous norm $\rho$;
$u\in{KS^{1,\alpha}(\Omega,\widetilde{G})}$ where $\widetilde{G}$
is another Carnot group with homogeneous norm $\widetilde{\rho}$;
$Du(p): G\rightarrow \widetilde{G}$ is the approximate Pansu
derivative of $u$ at $p\in{\Omega}$, see Definition \ref{pansu}
and Theorem \ref{pansuofSobolev}; $C$ is a constant and
$B_{c}(0,1)$ is the unit C-C ball centered at $0$. Our arguments
rely on the equivalence of several Sobolev classes between C-C
spaces. Let $R^{1,\alpha}(\Omega, M)$ and $H^{1,\alpha}(\Omega,M)$
denote the Sobolev spaces defined in the sense of Reshetnyak and
Haj\l asz respectively, see Definition \ref{sec2:definition1} and
Definition \ref{sec2:definition2}. When $\alpha>1$ and $\Omega$ is
a bounded open set in a C-C space with some conditions and $M$ is
a separable metric space, we prove that
\begin{equation}\label{intro:equivalence}
 KS^{1,\alpha}(\Omega,M)=R^{1,\alpha}(\Omega,M)=H^{1,\alpha}(\Omega,M)
\end{equation}
as sets, see Theorem \ref{sec2:equivalence}. The proof essentially
depends on several observations of the theory of real-valued
Sobolev classes defined on metric measure spaces which was
developed in \cite{KM}, \cite{FHK} and \cite{HK}. Let us mention
that  the equivalence of several definitions of Banach
space-valued Sobolev classes has been proven in \cite{HKST} where
an important technique, that each metric space $Y$ can be
isometrically embedded into a Banach space, for example into
$L^{\infty}(Y)$ or $l^{\infty}$ if $Y$ is separable, is trickily
adopted. Since such isometric embedding is not good enough (see
\cite{Semmes} for the fact that Heisenberg group is not
bilipschitz equivalent to any Euclidean space in any scale), we
will not use this idea. Compared with the proof suggested in
\cite{HKST}, our proof of (\ref{intro:equivalence}) is convenient
for our purpose, also direct and simpler due to the differential
structure of C-C spaces.

In \cite{Vo1} and \cite{Vo2}, Vodop'yanov made a systematic study
of $R^{1,\alpha}(\Omega, \widetilde{G})$ where $\alpha>1$,
$\Omega$ is a bounded open set of a Carnot group $G$ and
$\widetilde{G}$ is another Carnot group. In particular, he gave
several equivalent descriptions of
$R^{1,\alpha}(\Omega,\widetilde{G})$, including a characterization
using properties of coordinate functions which obviously covers
the first statement in Theorem \ref{intro:thmcl}. Equation
(\ref{intro:explicitcg}) is deduced from (\ref{intro:equivalence})
and the results in \cite{Vo1} and \cite{Vo2}. When $\Omega$ is an
Euclidean domain and $\widetilde{G}$ is the Heisenberg group,
(\ref{intro:explicitcg}) is just (\ref{intro:equcl}) (recall that
$R^{n}$ can be seen as an abelian Carnot group).

Although we can explicitly formulate the energy of Korevaar-Schoen
as (\ref{intro:explicitcg}), we do not know whether or not
$E^{\alpha}(\Omega,\widetilde{G})$ is lower-semicontinuous with
respect to some topology of $KS^{1,\alpha}(\Omega,\widetilde{G})$.
As done in \cite{KS}, the lower semicontinuity of
$E^{\alpha}(\Omega,M)$ with respect to the topology of
$L^{\alpha}(\Omega,M)$ is a byproduct of a subpartitional lemma
when $\Omega$ is a Riemannian domain, see also \cite{EF}. Sturm in
\cite{Sturm} generalized this fact to domains which possesses a
strong or weak ``measure contraction property", see also
\cite{KS1} and \cite{KS2}. Unfortunately, in general C-C spaces
seem to have no ``measure contraction property". We will
illustrate this fact for Heisenberg group in Section \ref{sec:3}.

Thus we will abandon the energy of Korevaar-Schoen. Instead we
will introduce the horizontal energy. Let us first recall the
definition of the energy in the theory of harmonic mappings
between smooth Riemannian manifolds (e.g. \cite{He},\cite{Jo3}).
Let $(M,g)$ and $(N,h)$ be two smooth manifolds with Riemannian
metric $g$ and $h$ respectively. The energy of a smooth map
$u:M\rightarrow N$ is defined as (up to a constant)
\begin{equation}
 E(u)=\int_{M}\|du\|^{\alpha}dv
\end{equation}
where $du$ is the induced differential map
$du(p):T_{p}M\rightarrow T_{f(p)}N$ ($du$ can be regarded as an
element of $\Gamma(T^{*}M\bigotimes u^{-1}TN)$); $\|du\|$ is the
norm with respect to the fiber metric of $\Gamma(T^{*}M\bigotimes
u^{-1}TN)$ induced by $u$ from $g$, $h$ and $dv$ is the volume
form in $M$. If we choose a coordinate chart of $M$ such that
$(\frac{\partial}{\partial x_{1}},\cdots,\frac{\partial}{\partial
x_{m} })$ is orthonormal with respect to $g$, then
(\ref{intro:riemainnianenergy}) can be rewritten as
\begin{equation}\label{intro:riemainnianenergy}
E(u)=\int_{M}\left (\sum_{i=1}^{m}h(du(\frac{\partial}{\partial
x_{i}}),du(\frac{\partial}{\partial x_{i}}))\right
)^{\frac{\alpha}{2}}dv.
\end{equation}

Now let $(G,\Delta,g_{c})$ and
$(\widetilde{G},\widetilde{\Delta},\widetilde{g}_{c})$ be two
sub-Riemannian manifolds. By definition, $G$ and $\widetilde{G}$
are two smooth manifolds endowed with smooth distributions
$\Delta=\textrm{span}\{X_{1},\cdots,X_{k}\}$,
$\widetilde{\Delta}=\textrm{span}\{Y_{1},\cdots,Y_{\widetilde{k}}\}$
respectively, and $g_{c}$ and $\widetilde{g}_{c}$ are fiberwise
inner products endowed to $\Delta$, $\widetilde{\Delta}$
respectively, such that $\{X_{1},\cdots,X_{k}\}$ and
$\{Y_{1},\cdots,Y_{\widetilde{k}}\}$ are orthonormal with respect
to $g_{c}$, $\widetilde{g}_{c}$ respectively. Note that any such
$g_{c}$ (or $\widetilde{g}_{c})$ can be realized as the
restriction of a Riemannian metric $g$ (or $\widetilde{g})$ on $G$
(or $\widetilde{G})$ to $\Delta$ (or $\widetilde{\Delta})$. Let
$u:G\rightarrow \widetilde{G}$ be a smooth map satisfying the
following contact condition
\begin{equation}\label{intro:contact}
  du(p)(X_{i}(p))\in{\widetilde{\Delta}_{(u(p))}}\quad \textrm{for }i=1,\cdots,m_{1}.
\end{equation}
We define the horizontal energy of $u$ as follows:
\begin{equation}\label{intro:contact-energy}
  HE(u)=\int_{G}\left
  (\sum_{i=1}^{k}\widetilde{g}(du(X_{i}),du(X_{i}))\right
  )^{\frac{\alpha}{2}}dv
\end{equation}
where $dv$ is the volume form in $G$ with respect to $\bar{g}$.

Note that $HE(u)$ is dependent on $g$ but independent of any
extension of $\widetilde{g}_{c}$. In the case $\Delta=TG$, $HE(u)$
only depends on $g_{c}$ and $\widetilde{g}_{c}$. The definition of
horizontal energy obviously generalizes the Riemannian energy
(\ref{intro:riemainnianenergy}) in the sense that if $\Delta=TG$
and $\widetilde{\Delta}=T\widetilde{G}$, then
(\ref{intro:contact-energy}) is just
(\ref{intro:riemainnianenergy}). Any smooth map satisfying
(\ref{intro:contact}) is called a contact map, see Definition
\ref{def:contact}. Any map in $R^{1,\alpha}(\Omega,\widetilde{G})$
satisfies (\ref{intro:contact}) in a weak sense, see Remark
\ref{remark:cha}. It turns out that
$R^{1,\alpha}(\Omega,\widetilde{G})$ is the natural space to study
the minimizing problem with respect to the horizontal energy. In
this paper, we will not explore the full general situation, but
restrict ourselves to C-C spaces, in particular to Carnot groups.
We will give an existence result of horizontal minimizers (see
Definition \ref{horizontalminimizer}) when the target is of Carnot
type.

In contrast to the easy existence problem of horizontal
minimizers, regularity problem is  very complicated. By now, we
have some results  in the case when $\Omega\subset R^{2}$ is
smooth and bounded open set and the target is the Heisenberg group
$H^{m}$. In this case, due to the conformal invariance of the
horizontal energy there is a close link to the two dimensional
isotropically constrained Plateau problem in $R^{2m}$ investigated
in \cite{SW} by Schoen-Wolfson when $m=2$ and in \cite{Qiu} by Qiu
Weiyang when $m>2$. The method of constructing isotropic
variations in \cite{Qiu} may be useful to further investigation.

To end this introduction, we sketch the structure of the paper. In
Section \ref{sec:1} we give notations, definitions and collect
some basic facts about C-C spaces and several definitions of
Sobolev classes defined in C-C spaces.  The equivalence of several
definitions of Sobolev classes from C-C spaces to separable metric
spaces will be proven in Section \ref{sec:2-1}, see Theorem
\ref{sec2:equivalence}. We discuss in \ref{sec:2-2} and
\ref{sec:2-3} the properties of $R^{1,\alpha}(\Omega,
\widetilde{G})$ such as several equivalent characterizations
(Theorem \ref{sec2:description}, \ref{sec2:characterization}),
precompactness (Theorem \ref{thm:campactness}) and the trace
problem (Theorem \ref{tracecharacterization}). In Section
\ref{sec:3} we discuss the properties of the energy of
Korevaar-Schoen (Theorem \ref{thm:energychara}) and give reasons
why we abandon it. We conjecture that C-C spaces do not possess
any type of ``measure contraction property". We will illustrate an
evidence to this conjecture by showing that Heisenberg group does
not possess the strong ``measure contraction property". So the
method used to prove that the approximate Korevaar-Schoen energies
satisfy a subpartitional lemma and then deduce that
Korevaar-Schoen energy is lower semicontinuous may not be valid in
this case. Section \ref{sec:4} is devoted to defining the
horizontal energy, to proving the existence of minimizers of the
horizontal energy minimizing problem when the domain is a smooth,
noncharacteristic bounded open set in a C-C space and the target
is a C-C space of Carnot type. The existence result is immediately
from the compactness theorem and the trivial lower semicontinuity
of the horizontal energy with respect to the weak topology.  In
Section \ref{sec:5} we discuss the regularity of the minimizers
when the domain is a bounded open set in $R^{2}$ and the target is
the Heisenberg group $H^{m}$.

{\small\noindent\textbf{Acknowledgement.} The author would like
very much to thank his advisor Professor Xiaoping Yang for his
constant encouragement and support in research as well as in
everyday life. He would also like to thank Professor S. K.
Vodop'yanov for sending the book ``Proceedings of analysis and
geometry, Novosibirsk: Sobolev Institute Press, 2000"  and several
his papers to him. He also thanks the referee for many valuable
comments and suggestions and reminding him of several important
papers about Sobolev functions in metric measure spaces.}

\section{Preliminaries and basic results}\label{sec:1}
The aim of this section is to fix the notations and collect some
basic results which will be used in the sequel.
\subsection{Carnot-Carath\'{e}odory spaces}\label{sec:1-1}
Let $\Delta =\textrm{span}\{X_{1},X_{2},\cdots,X_{k}\}$ be a
smooth distribution in $R^{n}$. We identify $X_{i}$ with a first
order differential operator in $R^{n}$.

Denote by $\overline{V}_{j}(p)$ the subspace of $T_{p}R^{n}=R^{n}$
spanned by all commutators of $X_{i}$'s of order $\leq j$
($\overline{V}_{1}=\Delta=\textrm{span}\{X_{1},\cdots,X_{k}\}$ is
called the \textit{horizontal bundle} whose cross sections are
called \textit{horizontal vector fields}). We say that $\Delta$
satisfies the \textit{H\"{o}rmander condition} provided for any
$p\in{R^{n}}$ there exists $r_{p}$ such that $\dim
(\overline{V}_{r_{p}}(p))=n$. $\Delta$ is \textit{equiregular} if
for each $j$, $\dim(\overline{V}_{j}(p))$ is independent of the
point. If $\Delta$ satisfies the H\"{o}rmander condition and is
equiregular, then the least integer $r$ such that
$\dim(\overline{V}_{r})=n$ is called the \textit{step} of
$\Delta$.

An absolutely continuous curve $\gamma:[a, b]\rightarrow R^{n}$ is
\textit{horizontal} if there exist Borel functions $c_{i}(t),a\leq
t\leq b$, such that
$\dot{\gamma}(t)=\sum_{i=1}^{k}c_{i}(t)X_{j}(\gamma(t))\textrm{
for a.e. }t\in{[a,b]}.$ We  endow a fiberwise inner product
$<\cdot,\cdot>_{c}$ to $\Delta$ such that $\{X_{1}(p),\cdots,
X_{k}(p)\}$ is orthonormal at every point $p\in{R^{n}}$. The
length of a horizontal curve $\gamma$ is defined as $
l_{c}(\gamma)=\int_{a}^{b}( \sum_{i=1}^{k}|c_{i}(t)|^{2})
^{\frac{1}{2}}dt $. Then the C-C distance $d_{c}$ between $p$ and
$q$ in $R^{n}$ is defined as the infimum of the lengths of all
horizontal curves connecting $p$ to $q$. $d_{c}$ is called the
\textit{C-C distance}. $R^{n}$ equipped with the C-C distance is
called \textit{C-C space}, denoted by $(R^{n},\Delta,d_{c})$. The
Chow theorem (\cite{Chow}) says that if the distribution $\Delta$
satisfies the H\"{o}rmander condition then there exists an
admissible curve connecting any given pair of points in $R^{n}$
and thus $d_{c}$ is a metric. For other equivalent definitions of
the C-C distance, we  refer to \cite{JSC}.
\begin{notations}\label{sec1:notation}
In the remainder of the paper, when we speak of a C-C space
$(R^{n},\Delta,d_{c})$ we assume that the distribution
$\Delta=\textrm{span}\{X_{1},\cdots,X_{k}\}$ satisfies the
H\"{o}rmander condition. We will use $\Delta_{p}$ to denote the
fiber of $\Delta$ through $p$. In the sequel $|E|$ will always
stand for $\mathcal{L}^{n}(E)$, where $\mathcal{L}^{n}$ is the
$n$-dimensional Lebesgue measure on $R^{n}$. $B_{c}(p,\delta)$
($B(p,\delta)$ or $B_{\delta}(p)$) will denote a C-C (Euclidean)
open ball centered at $p$ with radius $\delta$. We will use
$\overline{\Omega}$ to denote the closure of a subset
$\Omega\subset R^{n}$. By $\Omega\Subset\widetilde{\Omega}$ we
mean that $\overline{\Omega}$ is contained in
$\widetilde{\Omega}$. Let $u$ be a Borel function defined on
$\Omega\subset R^{n}$. The average value of $u$ on $\Omega$ will
be denoted by
$u_{\Omega}=\dashint_{\Omega}udx=|\Omega|^{-1}\int_{\Omega}udx$.
\end{notations}
\begin{lemma}[\cite{NSW}]
Let $(R^{n},\Delta,d_{c})$ be a C-C space. Then for every bounded
open set $\Omega\subset R^{n}$ there exists $C\geq 1$ such that
one has
\begin{equation}\label{sec1:doubling}
 |B_{c}(p,2\delta)|\leq C|B_{c}(p,\delta)|
\end{equation}
whenever $p\in{\Omega}$ and $\delta\leq 5\textrm{diam}\Omega$.
\end{lemma}
 The condition (\ref{sec1:doubling}) is
called the \textit{doubling condition} and the least constant $C$
such that (\ref{sec1:doubling}) holds is called the
\textit{doubling constant} and $Q:=\log_{2}C\geq n$ is called the
local \textit{homogeneous dimension} of $\Omega$. According to
\cite{Mi}, if $\Delta$ is equiregular, then the constant $
Q=\sum_{i=1}^{r}i(\dim(V^{i})-\dim(V^{i-1}))$ is the Hausdorff
dimension of $(R^{n},\Delta,d_{c})$. We refer to \cite{NSW} for
more about C-C balls.

Let $\Omega$ be a bounded open set in $(R^{n},\Delta,d_{c})$.
Following \cite{HeK} we say that a Borel function
$g:\Omega\rightarrow [0,\infty]$ is an \textit{upper gradient} of
another Borel function $u:\Omega\rightarrow R$ if for every
1-Lipschitz curve $\gamma:[0,T]\rightarrow \Omega$ we have
$
 |u(\gamma(0))-u(\gamma(T))|\leq \int_{0}^{T}g(\gamma(t))dt.
$ We recall that a curve $\gamma$ is called 1-Lipschitz if
$d_{c}(\gamma(t_{1}),\gamma(t_{2}))\leq |t_{2}-t_{1}|$ for all
$0\leq t_{1}<t_{2}\leq T$.

Let $u$ and $g\geq 0$ be two Borel functions defined on an open
subset $\Omega$. For the pair $(u,g)$ if there exist $C>0$ and
$\lambda\geq 1$ such that
\begin{equation}\label{sec1:1-pponcare}
  \dashint_{B_{c}}|u-u_{B_{c}}|dx\leq Cr\left(\dashint_{\lambda
B_{c}}g^{\alpha}\right)^{\frac{1}{\alpha}}
\end{equation}
holds for every metric ball $B_{c}$ in $\Omega$, where $r$ is the
radius of $B_{c}$, then we say the pair $(u,g)$ satisfies a
\textit{$(1,\alpha)$-Poincar\'{e} inequality for $C$ and
$\lambda$}.  We say $(R^{n},\Delta,d_{c})$ \textit{supports a
$(1,\alpha)$-Poincar\'{e} inequality}, $1\leq \alpha<\infty$, if
for every bounded open set $\Omega$ when $u$ is a continuous
function in $\Omega$ and $g$ is an upper gradient of $u$, the pair
$(u,g)$ satisfies a $(1, \alpha)$-Poincar\'{e} inequality for some
choice of constants $C_{\Omega}>0$ and $\lambda_{\Omega}\geq 1$.
The following theorem is well known, see \cite{Jerison},
\cite{JSC} and \cite{HK}.
\begin{theorem}\label{sec1:supports}
$(R^{n},\Delta,d_{c})$ supports a $(1,\beta)$-Poincar\'{e}
inequality for any $\beta\in{[1,\infty)}$.
\end{theorem}
For sharp results about Poincar\'{e} inequalities in metric
measure spaces we refer to \cite{Lu1}, \cite{Lu2} and \cite{FLW2}.
\begin{definition}[G-linear map]\label{glinearmap}
Let $G=(R^{n},V_{1},\delta_{\lambda})$ and
$\widetilde{G}=(R^{\widetilde{n}},\widetilde{V}_{1},\widetilde{\delta}_{\lambda})$
be two Carnot groups. A mapping $L:G\rightarrow \widetilde{G}$ is
called a G-linear map if
\begin{enumerate}
  \item L is a homogeneous with respect to $\delta_{\lambda}$ and
  $\widetilde{\delta}_{\lambda}$, that is,
  $L(\delta_{\lambda}p)=\widetilde{\delta}_{\lambda}L(p)$ for any $p\in{G}$ and $\lambda>0$.
  \item L is a group homomorphism, that is, $L(pq)=L(p)L(q)$ for
  any $p,q\in{G}$.
 \end{enumerate}
\end{definition}
Any G-linear map is smooth and contact, for a proof see e.g
\cite{Ma}.

\begin{definition}[Pansu differential]\label{pansu}
Let $G$ and $\widetilde{G}$ be two Carnot groups with homogeneous
norms $\rho$ and $\widetilde{\rho}$ respectively. Let $E$ be a
Borel subset of $G$. A G-linear map $L$ is called a Pansu
differential of a mapping $u:E\rightarrow \widetilde{G}$ at a
point $p\in{E}$ if
$$
 \lim_{x\rightarrow
p,x\in{E}}\frac{\widetilde{\rho}(L(a^{-1}x)^{-1}u(p)^{-1}u(x))}{\rho(p^{-1}x)}=0.
$$
 A G-linear map $L$
is called an approximate Pansu differential of $u$ in $E$ at a
point $p\in{U}$ if
$$
  \textrm{ap}\lim_{x\rightarrow p}\frac{\widetilde{\rho}
  (L(a^{-1}x)^{-1}u(p)^{-1}u(x))}{\rho(p^{-1}x)}=0
$$
where $\textrm{ap}\lim_{x\rightarrow p}f(x)$ denotes the
approximate limit of $f$ at $p$ (see \cite{Fe}).
\end{definition}
\begin{remark} The notion of derivatives for mappings between
Carnot groups was originally introduced by P. Pansu in
\cite{Pansu} where the set $E$ in Definition \ref{pansu} is
required to be an open set. The version of Definition \ref{pansu}
is due to \cite{VU} and \cite{Ma}.
\end{remark}

\subsection{Sobolev functions defined on Carnot-Carath\'{e}odory spaces}
There are several equivalent definitions for Sobolev functions on
metric measure spaces. The fundamental references in this topic
are \cite{HK}, \cite{Heinonen}. We concentrate on Sobolev
functions in C-C spaces. Due to the differential structure of C-C
spaces, the theory of Sobolev functions in C-C spaces are more
abundant than that in general metric measure spaces.

Let $(R^{n},\Delta,d_{c})$ be a C-C space and let $\Omega$ be an
open set in $R^{N}$. Let $\alpha$ be in $[1,\infty]$. The
\textit{horizontal Sobolev space} is the Banach space
$$
W_{X}^{1,\alpha}(\Omega)=\{u\in{L^{\alpha}(\Omega)|X_{i}u\in{L^{\alpha}(\Omega)}},j=1,\cdots,k\}
$$
endowed with the norm $\|
u\|_{W_{X}^{1,p}(\Omega)}=\|u\|_{L^{\alpha}(\Omega)}+\sum_{i=1}^{k}
\|X_{i}u\|_{L^{\alpha}(\Omega)}$. In the above definition,
$X_{i}u$ is understood in the distributional sense. Another way to
define the space $W_{X}^{1,\alpha}(\Omega)$ for
$1\leq\alpha<\infty$ is to take the closure of $C^{\infty}$
functions in the norm $\|\cdot\|_{W_{X}^{1,\alpha}(\Omega)}$. As
in the Euclidean case, the two approaches are equivalent. This was
obtained independently in  \cite{FSSC} and \cite{GN1}.

For $1\leq\alpha<\infty$, the Sobolev space $H^{1,\alpha}(\Omega)$
is defined as the set of all $u\in{L^{\alpha}(\Omega)}$ for which
there exists $0\leq g\in{L^{\alpha}(\Omega)}$ such that the
inequality
\begin{equation}\label{sec1:Hajlasz}
|u(x)-u(y)|\leq d_{c}(x,y)(g(x)+g(y))
\end{equation}
holds a.e. $x,y\in{\Omega}$. $H^{1,\alpha}(\Omega)$ is firstly
introduced by Haj\l asz in \cite{Ha}.  By $P^{1,\alpha}(\Omega)$
we denote the set of all functions $u\in{L^{\alpha}(\Omega)}$ such
that there exists $0\leq g\in{L^{\alpha}(\Omega)}$ such that the
pair $(u,g)$ satisfies a $(1,\alpha)$-Poincar\'{e} inequality.
Roughly speaking, the function $g$ in (\ref{sec1:Hajlasz})
corresponds to the maximal function of the gradient, while the
function $g$ in (\ref{sec1:1-pponcare}) looks more like the norm
of the gradient (see the Introduction of \cite{HK}). For other
notions of Sobolev functions in C-C spaces or general metric
measure spaces we refer the reader to \cite{Ch}, \cite{Sh},
\cite{FLW1}, \cite{LLW} and \cite{LP}.

The following theorem, which follows from Theorem
\ref{sec1:supports} and Theorem 1, Corollary 13 in \cite{FHK} (see
also \cite{KM}), is crucial to Theorem \ref{sec2:equivalence} in
Section \ref{sec:2}.
\begin{theorem}\label{sec1:equivalence}
Let $(R^{n},\Delta,d_{c})$ be a C-C space and $\Omega\subset
R^{n}$ be a bounded open set. Assume $1<\alpha<\infty$. Then the
following four conditions are equivalent.
\begin{enumerate}
  \item $u\in{W_{X}^{1,\alpha}(\Omega)}$.
  \item $u\in{H^{1,\alpha}(\Omega)}$.
  \item $u\in{P^{1,\alpha}(\Omega)}$.
  \item $u\in{L^{\alpha}(\Omega)}$ and there exist $0\leq
                  g\in{L^{\alpha}(\Omega)}$, constants $C>0$, $\lambda\geq 1$
                   such that $(u,g)$ satisfies a
                  $(1,\beta)$-Poincar\'{e} inequality for
                  $C, \lambda$ where $\beta\in{[1,\alpha)}$.
\end{enumerate}
Moreover
\begin{enumerate}
  \item[(i)]\, If $u\in{L^{\alpha}(\Omega)}$ and there exist $0\leq
                  g\in{L^{\alpha}(\Omega)}$, constants $C>0$, $\lambda\geq 1$
                  such that $(u,g)$ satisfies a
                  $(1,\beta)$-Poincar\'{e} inequality for $C$ and $\lambda\geq
                  1$ where $\beta\in{[1,\alpha)}$, then $
            \omega(x)=C(\sup_{r>0}\dashint_{B_{c}(x,r)}g^{\beta}(x)dx)^
{\frac{1}{\beta}}
$
is in $L^{\alpha}(\Omega)$ and the pair $(u,\omega)$ satisfies
(\ref{sec1:Hajlasz}) where $g$ is replaced by $\omega$.
  \item[(ii)]\, If $u\in{W_{X}^{1,\alpha}(\Omega)}$ and $(u,g)$
  satisfies a $(1,\alpha)$-Poincar\'{e} inequality, then $|Xu|\leq
  Cg$ a.e. for some constant $C$ independent of $u$ and $g$.
\end{enumerate}
\end{theorem}

\section{Sobolev classes from Carnot-Carath\'{e}odory spaces to separable metric
spaces}\label{sec:2} In this section we study Sobolev classes from
C-C spaces to separable metric spaces.  In Section \ref{sec:2-1}
we define $H^{1,\alpha},\,KS^{1,\alpha}$ and $R^{1,\alpha}$, then
we prove that they are equivalent as sets when $1<\alpha<\infty$.
In Section \ref{sec:2-2} we study properties of Sobolev mappings
from a C-C space to another C-C space of Carnot type by giving
several equivalent descriptions of $R^{1,\alpha}$ which slightly
generalizes some corresponding results in \cite{Vo1}.

\subsection{Equivalence of Sobolev classes}\label{sec:2-1}
Let $(R^{n},\Delta,d_{c})$ be a C-C space, $\Omega\subset R^{n}$
be a bounded open set with smooth boundary and $(M,d)$ be a
complete metric space with a (quasi-)metric $d$. Assume
$1\leq\alpha <\infty$. Let $u:\Omega\rightarrow M$ be a measurable
map. $u$ is called in $L^{\alpha}(\Omega,M)$ if $
\int_{\Omega}d(m_{0},u(p))^{\alpha}dp<\infty $ for some
$m_{0}\in{M}$. Since $\Omega$ is bounded, the definition is
independent of the choice of $m_{0}$ by the (quasi-)triangle
inequality of $d$. We identify two mappings which coincide
$\mathcal{L}^{n}$-almost everywhere. It is easily proved that
$L^{\alpha}(\Omega,M)$ is a complete metric space with the
distance $
d_{L^{\alpha}(\Omega,M)}(u,v)=\int_{\Omega}d(u(p),v(p))^{\alpha}dp,
$ see e.g. \cite{RM}.

For $\epsilon>0$, let $
\Omega_{\epsilon}:=\{p\in{\Omega}:\textrm{dist}_{c}(p,\partial\Omega)>\epsilon\}
$ with
$$
\textrm{dist}_{c}(p,\partial\Omega)=\inf_{q\in{\partial\Omega}}{d_{c}(p,q)}.
$$
For a map $u:\Omega\rightarrow M$ and for a point $p\in{\Omega}$,
we define the averaged $\epsilon$-approximate density function
$$
e_{\epsilon}^{\alpha}(p;u)=\dashint_{B_{c}(p,\epsilon)}\left(\frac{d(u(p),u(q))}{\epsilon}\right)^{\alpha}dq
$$
 where $1\leq\alpha<\infty$ and $p\in{\Omega}$. If
$\varphi\in{C_{c}(\Omega,[0,1])}$ and
$\epsilon<\textrm{dist}_{c}(\textrm{supp}\varphi,\partial\Omega)$,
we define the approximate energy
\begin{equation}\label{approximateenergy}
E_{\epsilon}^{\alpha}(\varphi;u)=\int_{\Omega}\varphi(p)e_{\epsilon}^{\alpha}(p;u)dp.
\end{equation}
 We now define the class $KS^{1,\alpha}(\Omega,M)$.
\begin{definition}
Let $u\in{L^{\alpha}(\Omega,M)}$. We say $u$ is in
$KS^{1,\alpha}(\Omega,M)$ if
$$
E^{\alpha}(u,\Omega)=\sup_{\varphi\in{C_{c}(\Omega,[0,1])}}
\limsup_{\epsilon\rightarrow 0}E_{\epsilon}^{\alpha}(\varphi;u)
$$
is finite. If $u\in{KS^{1,\alpha}(\Omega,M)}$,
$E^{\alpha}(u,\Omega)$ is called the \textit{energy of
Korevaar-Schoen}.
\end{definition}
The above definition is firstly introduced in \cite{KS} by
Korevaar and Schoen in the case where $\Omega$ is a Riemannian
domain. Later it is generalized to general metric measure spaces,
see \cite{RM} and \cite{HKST}.

\begin{definition}\label{sec2:definition1} Let $u\in{L^{\alpha}(\Omega,M)}$. We say
$u\in{H^{1,\alpha}(\Omega,M)}$ if there exists
$0\leq\omega\in{L^{\alpha}(\Omega)}$ such that
\begin{equation}\label{eq:hajlasz}
d(u(p),u(q))\leq d_{c}(p,q)(\omega(p)+\omega(q))
\end{equation}
holds for a.e. $p,q\in{\Omega}$. We set
$$
E^{\alpha}_{H}(u,\Omega)=\inf_{\omega}\|\omega\|^{\alpha}_{L^{\alpha}(\Omega)}
$$
where the infimum is taken among all nonnegative functions
$\omega$ in $L^{\alpha}(\Omega)$ such that \eqref{eq:hajlasz}
holds.
\end{definition}
$H^{1,\alpha}(\Omega,M)$ is a natural generalization of
$H^{1,\alpha}(\Omega)$ and can also be extended to more general
metric measure spaces (\cite{HKST}).

\begin{definition}\label{sec2:definition2} Let $u\in{L^{\alpha}(\Omega,M)}$. We say
$u\in{R^{1,\alpha}(\Omega,M)}$ if for any $m\in{M}$, the scalar
function $\theta_{m}(p)$ defined by $\theta_{m}(p):=d(m,u(p))$ is
in $W_{X}^{1,\alpha}(\Omega)$ and there exists $0\leq
g\in{L^{\alpha}(\Omega)}$ (independent of $m$) such that
$$
|X\theta_{m}(p)|\leq g(p)
$$
a.e. $p\in{\Omega}$ for any $m\in{M}$. We call $g$ is a dominant
function of $u$. We set
$$
E^{\alpha}_{R}(u,\Omega)=\inf_{g}\|g\|^{\alpha}_{L^{\alpha}(\Omega)}
$$
where the infimum is taken among all dominant functions $g$ of
$u$.
\end{definition}
When $\Omega$ is an Euclidean domain (that is
$\Delta=\textrm{span}\{\frac{\partial}{\partial
x_{1}},\cdots,\frac{\partial}{\partial x_{n}}\}$), this definition
coincides with that in \cite{Re}.

\begin{lemma}\label{sec2:poincare} If
$u\in{KS^{1,\alpha}(\Omega,M)}$, then
$\theta_{m}\in{KS^{1,\alpha}(\Omega,R)}$ and $
e_{\epsilon}^{\alpha}(p,\theta_{m})\leq
Ce_{\epsilon}^{\alpha}(p,u)$ for any $m\in{M}$ and some constant
$C$ (independent of $m$). Thus there exists $0\leq
g\in{L^{\alpha}(\Omega)}$ independent of $m$ such that the pair
$(\theta_{m},g)$ satisfies a $(1,\beta)$-Poincar\'{e} inequality
for any $m\in{M}$ and $\beta\in{[1,\alpha)}$.
\end{lemma}
Lemma \ref{sec2:poincare} is from Theorem \ref{sec1:1-pponcare}
and a careful examination of the proof of Theorem 4.5 in
\cite{KM}.

Various seemly different definitions of Sobolev classes are
equivalent. Precisely, we have
\begin{theorem}\label{sec2:equivalence}
Let $(R^{n},\Delta,d_{c})$ be a C-C space and $\Omega$ be a
bounded open set in $R^{n}$. Assume that $M$ is a separable metric
space with a (quasi-) metric $d$. If $1<\alpha<\infty$, then
$$
KS^{1,\alpha}(\Omega,M)=R^{1,\alpha}(\Omega,M)=H^{1,\alpha}(\Omega,M)
$$
 as sets.
\end{theorem}
\begin{proof}
 \textbf{ Step 1.} $KS^{1,\alpha}(\Omega,M)\subset
R^{1,\alpha}(\Omega, M)$.

Let $u\in{KS^{1,\alpha}(\Omega,M)}$, then from Theorem
\ref{sec1:equivalence} and Lemma \ref{sec2:poincare} we conclude
that $\theta_{m}\in{W_{X}^{1,\alpha}(\Omega)}$ and there exists
$g^{\prime}\in{L^{\alpha}(\Omega)}$ independent of $m$ such that $
|X\theta_{m}(p)|\leq g^{\prime}(p)$ for any $m\in{M}$ and a.e.
$p\in{\Omega}$. Thus $KS^{1,\alpha}(\Omega,M)\subset
R^{1,\alpha}(\Omega, M)$.

 \textbf{Step 2.} $R^{1,\alpha}(\Omega,M)\subset H^{1,\alpha}(\Omega,M)$.

Let $u\in{R^{1,\alpha}(\Omega,M)}$. By definition,
$u\in{L^{\alpha}(\Omega,M)},\,\theta_{m}\in{W_{X}^{1,\alpha}(\Omega)}$
for any $m\in{M}$ and there exists $g\in{L^{\alpha}(\Omega)}$
independent of $m$ such that $|X\theta_{m}(p)|\leq g(p)$ holds
a.e. $p\in{\Omega}$ for any $m\in{M}$. Since the pair
$(\theta_{m},X\theta_{m})$ satisfies a $(1,\alpha)$-Poincar\'{e}
inequality for some choice of constants $C>0$ and $\lambda\geq 1$
for any $m\in{M}$ (see e.g. \cite{GN}), the pair $(\theta_{m},g)$
also satisfies a $(1,\alpha)$-Poincar\'{e} inequality for $C$ and
$\lambda$ and any $m\in{M}$. For $\beta\in{[1,\alpha)}$ we let
$$
\omega(x)=C\left(\sup_{r>0}\dashint_{B_{c}(x,r)}g^{\beta}(x)dx\right)^
{\frac{1}{\beta}}.
$$
 Since $M$ is separable, we can choose a
sequence of points $\{m_{i}\}_{i=1}^{\infty}$ such that
$\{m_{i}\}_{i=1}^{\infty}$ is dense in $M$. From Theorem
\ref{sec1:equivalence} we conclude that for any $i$ there exists a
set $\Omega_{i}$ of measure zero such that the inequality
$$
  |\theta_{m_{i}}(p)-\theta_{m_{i}}(q)|\leq d_{c}(p,q)( \omega(p)+\omega(q))
$$
 holds for any $p,q\in{\Omega\backslash\Omega_{i}}$. Let us set
$\Omega^{\prime}:=\bigcup_{i=1}^{\infty}\Omega_{i}$, then
$|\Omega^{\prime}|=0$ and
\begin{equation}\label{sec2:hajlasz}
 |d(m_{i},u(p))-d(m_{i},u(q))|\leq d_{c}(p,q)( \omega(p)+\omega(q))
\end{equation}
holds for any $p,q\in{\Omega\backslash \Omega^{\prime}}$ and
$i\in{\mathbb{N}}$. Fixing a point $p$ such that
(\ref{sec2:hajlasz}) holds for any $q\in{\Omega\backslash
\Omega^{\prime}}$ and $i\in\mathbb{N}$, we can choose a
subsequence of $\{m_{i}\}_{i=1}^{\infty}$ such that it converges
to $u(p)$. Hence $ d(u(p),u(q))\leq d_{c}(p,q)(
\omega(p)+\omega(q))$ holds for all $p,q\in{\Omega\backslash
\Omega^{\prime}}$. Thus $u\in{H^{1,\alpha}(\Omega,M)}$. \vskip
10pt

 \textbf{Step 3.} $H^{1,\alpha}(\Omega,M)\subset KS^{1,\alpha}(\Omega,M)$.

Let $u\in{H^{1,\alpha}(\Omega,M)}$. By definition,
$u\in{L^{\alpha}(\Omega,M)}$ and there exists
$0\leq\omega\in{L^{\alpha}(\Omega)}$ such that
\begin{equation}\label{sec2:Hajlasz}
 d(u(p),u(q))\leq d_{c}(p,q)(\omega(p)+\omega(q))
\end{equation}
holds for a.e. $p,q\in{\Omega}$.  Let
$\varphi\in{C_{c}(\Omega,[0,1])}$ and
$\epsilon<\frac{1}{3}\textrm{dist}_{c}(\textrm{supp}\varphi,\partial\Omega)$.
We have
  \begin{align}
  E_{\epsilon}^{\alpha}(\varphi;u) & \leq \int_{\Omega}\varphi
  \dashint_{B_{c}(p,\epsilon)}\left(\frac{d(u(p),u(q))}{\epsilon}\right)^{\alpha}dqdp\label{0}\\
  \quad & \leq \int_{\Omega}\varphi
  \dashint_{B_{c}(p,\epsilon)}\left(\frac{d(u(p),u(q))}{d_{c}(p,q)}\right)^{\alpha}dqdp\notag\\
 \quad &\leq
  C\int_{\Omega}\varphi\dashint_{B_{c}(p,\epsilon)}\left(|\omega(p)|^{\alpha}+
  |\omega(q)|^{\alpha}\right)dqdp \label{1}\\
  \quad & \leq
  C\|\omega\|^{\alpha}_{L^{\alpha}(\Omega)}+C\int_{\Omega_{3\epsilon}}\left(\dashint_{B_{c}(p,\epsilon)}|
  \omega(q)|^{\alpha}dq\right)dp\notag\\
  \quad & =C\|\omega\|^{\alpha}_{L^{\alpha}(\Omega)}+C\int_{\Omega_{3\epsilon}}\left(\int_{\Omega_{2\epsilon}}|
  \omega(q)|^{\alpha}\frac{\chi_{B_{c}(p,\epsilon)}(q)}{|B_{c}(p,\epsilon)|}dq\right)dp\notag\\
  \quad &\leq C\|\omega\|^{\alpha}_{L^{\alpha}(\Omega)}+C\int_{\Omega}\left(\int_{\Omega_{2\epsilon}}|
  \omega(q)|^{\alpha}\frac{\chi_{B_{c}(q,\epsilon)}(p)}{|B_{c}(p,\epsilon)|}dq\right)dp \label{2}\\
  \quad &= C\|\omega\|^{\alpha}_{L^{\alpha}(\Omega)}+C\int_{\Omega_{2\epsilon}} |\omega(q)|^{\alpha}
  \left(\int_{\Omega}
 \frac{\chi_{B_{c}(q,\epsilon)}(p)}{|B_{c}(p,\epsilon)|}dp\right)dq \label{3}\\
 \quad &= C\|\omega\|^{\alpha}_{L^{\alpha}(\Omega)}+C\int_{\Omega_{2\epsilon}} |\omega(q)|^{\alpha}\left(\int_{B_{c}(q,\epsilon)}
 \frac{1}{|B_{c}(p,\epsilon)|}dp\right)dq\notag\\
 \quad &= C\|\omega\|^{\alpha}_{L^{\alpha}(\Omega)}+C\int_{\Omega_{2\epsilon}} |\omega(q)|^{\alpha}
 \left(\dashint_{B_{c}(q,\epsilon)}
 \frac{|B_{c}(q,\epsilon)|}{|B_{c}(p,\epsilon)|}dp\right)dq\notag\\
 \quad &\leq C\|\omega\|^{\alpha}_{L^{\alpha}(\Omega)}+C\int_{\Omega_{2\epsilon}} |\omega(q)|^{\alpha}
 \left(\dashint_{B_{c}(q,\epsilon)}
 \frac{|B_{c}(p,2\epsilon)|}{|B_{c}(p,\epsilon)|}dp\right)dq \label{4}\\
 \quad &\leq C^{\prime}\|\omega\|^{\alpha}_{L^{\alpha}(\Omega)}\label{5}
  \end{align}
where in \eqref{1} we used \eqref{sec2:Hajlasz}; \eqref{2} is from
$\chi_{B_{c}(p,\epsilon)}(q)=\chi_{B_{c}(q,\epsilon)}(p)$ where
$\chi_{A}(q)$ denotes the characteristic function of the set $A$;
in \eqref{3} we used the Fubini's Theorem; \eqref{4} is from the
fact that if $p\in{B_{c}(q,\epsilon)}$ then
$B_{c}(q,\epsilon)\subset B_{c}(p,2\epsilon)$; \eqref{5} is from
the doubling condition \eqref{sec1:doubling}.

So $u\in{KS^{1,\alpha}(\Omega,M)}$.
\end{proof}
\begin{corollary}\label{energyequivalence}
Let $(R^{n},\Delta,d_{c})$ be a C-C space and $\Omega$ be a
bounded open set in $R^{n}$. Assume that $M$ is a separable metric
space with a (quasi-) metric $d$. If $1<\alpha<\infty$ and
$u\in{R^{1,\alpha}(\Omega,M)}$, then
$E^{\alpha}(u,\Omega),E_{H}^{\alpha}(u,\Omega)$ and
$E_{R}^{\alpha}(u,\Omega)$ are equivalent in the sense that each
one can be dominated by a constant multiple of another.
\end{corollary}

\subsection{Basic properties of Sobolev mappings}\label{sec:2-2}
In this section we slightly generalize the results in \cite{Vo1}
of some equivalent descriptions of
$R^{1,\alpha}(\Omega,\widetilde{G})$ where $\Omega\subset G$ is a
bounded open set and $G$, $\widetilde{G}$ are two Carnot groups to
the case when $G$ is a C-C space and $\widetilde{G}$ is a C-C
space of Carnot type (see Definition \ref{sec2:carnottype}).

\textbf{In the sequel we will assume that $\Delta$ is
equiregular}.
\begin{definition}\label{sec2:carnottype}
A C-C space $(R^{n},\Delta,d_{c})$ is \textit{of Carnot type} if
the system $\Delta=\textrm{span}\{X_{1},\cdots,X_{k}\}$ is of the
form
\begin{equation}\label{eq:carnottype}
  X_{i}(p)=\frac{\partial}{\partial
x_{i}}+\sum_{j=k+1}^{n}a_{i}^{j}(p)\frac{\partial}{\partial x_{j}}
\quad i=1,\cdots,k,
\end{equation}
where $a_{i}^{j}$ are smooth.
\end{definition}
This definition is motivated by the analogy with the canonical
generating vector fields of a Carnot group (see
\eqref{sec1:coefficients}).

 Next we will use the concept
of \textit{``some property holds for a.e. curves"}. Let us briefly
describe it, for details see \cite{KR}. Let $(R^{n},\Delta,d_{c})$
be a C-C space and $A\subset R^{n}$ be a bounded open set. Let
$\Gamma$ be a fibration of $A$ satisfying that the role of a fiber
$\gamma\in{\Gamma}$ is played by integral curves of a vector field
$\tau\in{\textrm{span}\{X_{1},\cdots,X_{k}\}}$. If we denote the
flow induced by the field by the symbol $f_{s}$ then the fiber has
the form $\gamma(s)=f_{s}(p)$, where $p$ belongs to a hypersurface
$\Sigma$ transversal to $\tau$ (such $\Sigma$ exists obviously).
We can endow a measure $d\gamma$ to $\Gamma$ as follows
$$d\gamma=\mathcal{F}_{f_{-s}}i(\tau)dx$$
where $\mathcal{F}_{f_{s}}$ is the Jacobian of the flow $f_{s}$
,$i(\tau)$ is the interior product of the vector field $\tau$ and
$dx$ is the standard volume form in $R^{n}$, such that
$$
c_{0}|B|^{\frac{Q-1}{Q}}\leq \int_{\gamma\in{\Gamma},\gamma\cap
B_{c}(x,r)\neq \varnothing}d\gamma\leq c_{1}|B|^{\frac{Q-1}{Q}}.
$$
for sufficiently small balls $B=B_{c}(x,r)\subset R^{n}$ with
constants $c_{0}$ and $c_{1}$, where $Q$ is the homogeneous
dimension of $(R^{n},\Delta,d_{c})$. We can identify a fiber of
$\Gamma$ with  a point in $\Sigma$ through the canonical
projection. Roughly speaking, saying that some property holds for
a.e. curves in $\Gamma$ is the same as saying that this property
holds for $d\sigma$ a.e. points in $\Sigma$ where $d\sigma$
denotes the Riemannian measure on $\Sigma$ induced from the
standard Euclidean metric in $R^{n}$. \vskip10pt

\begin{definition}[$ACL(\Omega,M)$]\label{sec2:ACL}
Let $(R^{n},\Delta,d_{c})$ be a C-C space and $\Omega\subset
R^{n}$ be a bounded open set. Let $M$ be a metric space with a
(quasi-)metric $d$. A mapping $u:\Omega\rightarrow M$ is
\textit{absolutely continuous on lines} (denoted by ACL for
brevity) if for every fibration $\Gamma_{i}$ of $\Omega$
 determined by $X_{i}$, $i=1,\cdots,k$, the curve
$u(\gamma): \gamma\cap\Omega\rightarrow M$, is absolutely
continuous in the parameter $t$ for $d\gamma$-almost every curve
$\gamma\in{\Gamma_{i}}$.
\end{definition}

In Definition \ref{sec2:ACL} any element $\gamma$ in $\Gamma_{i}$
is a flow induced by $X_{i}$. Since $\Omega$ is bounded,
$\gamma\cap \Omega$ has the form $\exp_{p}(tX_{i})$ where
$p\in{\Omega}$ and vice-versa.
\begin{definition}[\textbf{$R_{1}^{1,\alpha}(\Omega,M)$ and
$R_{2}^{1,\alpha}(\Omega,M)$}]\label{r1} Let
$(R^{n},\Delta,d_{c})$ be a C-C space and $\Omega\subset R^{n}$ be
a bounded open set. Let $M$ be a separable metric space with a
(quasi-)metric $d$. Assume $u:\Omega\rightarrow M $ be a mapping
and $1\leq\alpha<\infty$. We say
$u\in{R_{1}^{1,\alpha}(\Omega,M)}$ if
\begin{enumerate}
\item $\theta_{m}\in{L^{\alpha}(\Omega)}$ for any $m\in{M}$.
\item up to a modification on a set of measure zero, $u\in{ACL(\Omega,M)}$;
moreover the length of the curve $u(\gamma):\gamma\cap\Omega\rightarrow
  M$ is absolutely continuous in the parameter $t$ for $d\gamma$
  a.e. curve $\gamma\in{\Gamma_{i}}$ where $\Gamma_{i}$ is a
  fibration of $\Omega$ determined by $X_{i}$.
\item the derivative $
  X_{i}u_{l}(p)=\lim_{t\rightarrow 0}\frac{l(u(p),u(\exp_{p}(tX_{i})))}{t}$
  of the length of the curve $\Upsilon(\tau):=u(\exp_{p}(\tau X_{i})):[0,t]\rightarrow
  M,$ which exists almost everywhere in $\Omega$, belongs to
  $L^{\alpha}(\Omega)$ for all $i=1,\cdots,k$. Here
  $l(u(p),u(\exp_{p}(tX_{i})))$ is the length of the path
  $\Upsilon[0,t]$.
\end{enumerate}
We say $u\in{R_{2}^{1,\alpha}(\Omega,M)}$ if
\begin{enumerate}
\item for any function $f\in\textrm{Lip}(M)$, $f\circ u\in{W_{X}^{1,\alpha}
  (\Omega)}$;
\item there exists $0\leq g\in{L^{\alpha}(\Omega)}$ such that
  $|X(f\circ u)|\leq \textrm{Lip}f. g$ holds a.e. for all
  $f\in{\textrm{Lip}(\Omega)}$.
 \end{enumerate}
\end{definition}

The following theorem is a slight generalization of Proposition
4.1 in \cite{Vo1} (see also \cite{VU}).

\begin{theorem}\label{sec2:description}
Let $(R^{n},\Delta,d_{c})$ be a C-C space and $\Omega\subset
R^{n}$ be a bounded open set. Let $M$ be a separable metric space
with a (quasi-)metric $d$. Assume $u:\Omega\rightarrow M $ be a
mapping and $1\leq\alpha<\infty$, then
$$
R^{1,\alpha}(\Omega,M)=R_{1}^{1,\alpha}(\Omega,M)=R_{2}^{1,\alpha}(\Omega,M)
$$
as sets.
\end{theorem}

\begin{definition}[\textbf{$HW^{1,\alpha}(\Omega,M)$ and
$HW_{1}^{1,\alpha}(\Omega,M)$}]\label{w1} Let
$(R^{n},\Delta,d_{c})$ be a C-C space and $\Omega\subset R^{n}$ be
a bounded open set. Let
$M=(R^{\widetilde{n}},\widetilde{\Delta},\widetilde{d})$ be a C-C
space with
$\Delta=\textrm{span}\{\widetilde{X}_{1},\cdots,\widetilde{X}_{\widetilde{k}}\}$
and a (quasi-)metric $\widetilde{d}$ which is equivalent to the
C-C metric $\widetilde{d}_{c}$ (that is, there exist constants
$C_{1}$ and $C_{2}$ such that $C_{1}\widetilde{d}\leq
\widetilde{d}_{c}\leq C_{2}\widetilde{d}$). Assume
$u=(u^{1},\cdots,u^{\widetilde{n}}):\Omega\rightarrow M$ be a
mapping and $1\leq\alpha<\infty$. We say
$u\in{HW^{1,\alpha}(\Omega,M)}$ if
\begin{enumerate}
  \item $\widetilde{d}(u)\in{L^{\alpha}(\Omega)}$;
  \item up to a redefinition on a set of measure zero, $u^{i}\in{ACL(\Omega,R)}$ for
  $i=1,\cdots,\widetilde{n}$;
  for $i=1,\cdots,\widetilde{k}$;
  \item
  $
         X_{j}u(x)=\sum_{i=1}^{\widetilde{n}}X_{j}u^{i}(x)\frac{\partial}{\partial
         \widetilde{x}_{i}}\in{\widetilde{\Delta}_{u(x)}}
 $
  which exists for a.e
  $x\in{\Omega}$, belongs to $L^{\alpha}(\Omega)$, that is,
  $\int_{\Omega}\|X_{j}u(x)\|_{\widetilde{<\cdot,\cdot>}_{\widetilde{c}}}^{\alpha}dx<\infty$,
  $j=1,\cdots,
  \widetilde{k}$, where $\widetilde{<\cdot,\cdot>}_{\widetilde{c}}$ denotes the fiberwise inner product
  in $\widetilde{\Delta}$.
  \end{enumerate}
   We say $u\in{HW_{1}^{1,\alpha}(\Omega,M)}$ if
  \begin{enumerate}
  \item $\widetilde{d}(u)\in{L^{\alpha}(\Omega)};$
  \item up to a redefinition on a set of measure zero, $u\in{ACL(\Omega,M)}$;
  \item the derivative
 $
  X_{j}u(x)=\frac{d}{dt}u(\exp_{x}(tX_{j}))\mid_{t=0}\in{\widetilde{\Delta}_{u(x)}}
$
 which exists a.e. in $\Omega$, belongs to $L^{\alpha}(\Omega)$, $j=1,\cdots,k$.
\end{enumerate}
\end{definition}
In Definition \ref{w1}, we abuse the notation
$\widetilde{d}(\widetilde{x})=\widetilde{d}(\widetilde{x},0)$ for
$\widetilde{x}\in{R^{\widetilde{n}}}$. \vskip 10pt

The following lemma, draw from \cite{Hale},  is crucial to prove
Theorem \ref{sec2:characterization}.
\begin{lemma}[Carath\'{e}odory]\label{lemma:caratheodory}
Suppose $D$ is an open set in $R^{N+1},$ $f(t,x):D\rightarrow
R^{N}$ satisfies the Carath\'{e}odory conditions on $D$, that is,
$f$ is Borel measurable in $t$ and for each compact set
$D^{\prime}$ of $D$, there is an integrable function
$m_{D^{\prime}}$ such that $ |f(t,x)|\leq
m_{D^{\prime}}(t),(t,x)\in{D^{\prime}}.$ Moreover for each compact
set $U$ in $D$, there exists an integrable function $k_{U}(t)$
such that $ |f(t,x)-f(t,y)|\leq k_{U}(t)|x-y|, (t,x)\in{U},
(t,y)\in{U}.$ Then, for any $(t_{0},x_{0})$ in $U$, there exists a
unique solution $x(t,t_{0},x_{0})$ of $ \frac{dx(t)}{dt}=f(t,x)
\textrm{ a.e }$ passing trough $(t_{0},x_{0})$. Moreover the
domain $E$ in $R^{N+2}$ of the function $x(t,t_{0},x_{0})$ is open
and $x(t,t_{0},x_{0})$ is continuous in $E$.
\end{lemma}

The following theorem, which can be proved by using Lemma
\ref{lemma:caratheodory} and a similar argument of S. K.
Vodop'yanov (Proposition 4.2 in \cite{Vo1}), is of paramount
importance for our purpose.

\begin{theorem}\label{sec2:characterization}
Let $(R^{n},\Delta,d_{c})$ be a C-C space and $\Omega\subset
R^{n}$ be a bounded open set. Let
$M=(R^{\widetilde{n}},\widetilde{\Delta},\widetilde{d})$ be a C-C
space of Carnot type where
$\Delta=\textrm{span}\{\widetilde{X}_{1},\cdots,\widetilde{X}_{\widetilde{k}}\}$
 and $\widetilde{d}$ is a metric equivalent to the
C-C metric $\widetilde{d}_{c}$. Assume
$u=(u^{1},\cdots,u^{\widetilde{n}}):\Omega\rightarrow M$ be a
mapping and $1\leq\alpha<\infty$. Then
$$
R^{1,\alpha}(\Omega,M)=R_{1}^{1,\alpha}(\Omega,M)=
R_{2}^{1,\alpha}(\Omega,M)=HW^{1,\alpha}(\Omega,M)=
HW_{1}^{1,\alpha}(\Omega,M)
$$
as sets.
\end{theorem}

\begin{remark}\label{generalccspace}
 Note that if $u(\exp_{p}(tX_{j}))$ is a horizontal
curve in
$M=(R^{\widetilde{n}},\widetilde{\Delta},\widetilde{d}_{c})$, then
it follows from
$$
l_{c}(u(\exp_{p}(t_{1}X_{j})),u(\exp_{p}(tX_{j})))=\int_{t_{1}}^{t}\|X_{j}u(\exp_{p}(sX_{j}))\|_{\widetilde{<\cdot,\cdot>}_{c}}ds
$$
that $
X_{j}u_{l_{c}}(\exp_{p}(t_{1}X_{j})=\|X_{j}u(\exp_{p}(t_{1}X_{j}))\|_{\widetilde{<\cdot,\cdot>}_{c}}
\textrm{ a.e. }t_{1},$ where the length $l_{c}$ of
$u(\exp_{p}(tX_{j}))$ is computed with respect to the C-C metric
$\widetilde{d}_{c}$. So if $u(\exp_{p}(tX_{j}))$ is horizontal,
then $C_{1}X_{j}u_{l}(p)\leq
\|X_{j}u(p)\|_{\widetilde{<\cdot,\cdot>}_{c}}\leq
C_{2}X_{j}u_{l}(p)$ trivially holds for a.e. $p\in{\Omega}$ where
$C_{1},C_{2}$ are constants only depends on the (quasi)-metric
$\widetilde{d}$.
\end{remark}
\begin{definition}[contact mapping]\label{def:contact}
Let $(R^{n},\Delta,d_{c})$ and
$(R^{\widetilde{n}},\widetilde{\Delta},\widetilde{d}_{c})$ be two
C-C spaces. Let $\Omega$ be a bounded open set of $R^{n}$. Assume
$u:\Omega\rightarrow R^{\widetilde{n}}$ be a measurable mapping.
We say $u$ is a weakly contact mapping if
\begin{enumerate}
  \item
  $X_{j}u^{i}(p)=\frac{du^{i}(\exp_{p}{tX_{j}})}{dt}\mid_{t=0}$
  exists a.e. $p\in{\Omega}$ for $j=1,\cdots,k$ and
  $i=1,\cdots,\widetilde{n}$;
  \item\label{contact} $X_{j}u(p)=\sum_{i=1}^{\widetilde{n}}X_{j}u^{i}(p)\frac{\partial}
  {\partial \widetilde{x}_{i}}\in{\widetilde{\Delta}_{u(p)}}$ a.e.$p\in{\Omega}$.
 \end{enumerate}
If $u$ is smooth and satisfies (\ref{contact}), then $u$ is called
a contact map. If $u$ is a (weakly) contact mapping, then $u$
induces a linear map $ D_{h}u(p):\Delta_{p}\rightarrow
\widetilde{\Delta}_{u(p)}. $
\end{definition}

\begin{remark}\label{remark:cha}Under the same conditions of Theorem \ref{sec2:characterization}
two observations are in order:
\begin{enumerate}
\item If $u\in{R^{1,\alpha}(\Omega,M)}$, then by Theorem
  \ref{sec2:characterization}, $u$ is a weakly contact mapping and
  the induced map $D_{h}u$ can be represented by the matrix
  $(X_{i}u^{j}(p))_{k\times \widetilde{k}}$ of which each entry
  belongs to $L^{\alpha}(\Omega)$. It is easily inferred from \eqref{eq:carnottype} that $
  X_{j}u(p)=\sum_{i=1}^{\widetilde{k}}X_{j}u^{i}(p)\widetilde{X}_{i}.$
\item\label{item3} Let $u\in{R^{1,\alpha}(\Omega,M)}$. By
 Theorem \ref{sec2:characterization} we have
 $u\in{L^{\alpha}(\Omega,M)}$,
$u^{i}\in{ACL(\Omega,R)}$ for $i=1,\cdots,\widetilde{n}$ and
$X_{j}u^{i}\in{L^{\alpha}(\Omega)}$ for
$i=1,\cdots,\widetilde{k}$. We can not verify
\begin{equation}\label{integrable2}
  u^{i}\in{L^{\alpha}(\Omega)} \quad\textrm{for } i=1,\cdots,\widetilde{k}.
\end{equation}
 But if
$M$ is a Carnot group with a homogeneous norm $\widetilde{\rho}$,
then \eqref{integrable2} holds. In general we do not have  that
$u^{i}\in{W_{X}^{1,\alpha}(\Omega)}$ for
$i=\widetilde{k}+1,\cdots,\widetilde{n}$ even if $M$ is a Carnot
group.
  \end{enumerate}
\end{remark}

In section \ref{sec:3} we will use the results about Pansu
differentiability of Sobolev mappings between Carnot groups with
respect to the topology of $L^{\alpha}(\Omega)$ to get the
explicit form of the Korevaar-Schoen energy. Pansu
differentiability with respect to several topology for Sobolev
mappings between Carnot groups has been studied in details in
\cite{Vo1}, \cite{Vo2}, \cite{Vo3} and \cite{VU}.

\begin{theorem}\label{pansuofSobolev}
Let $G=(R^{n},V_{1},\delta_{\lambda},\rho)$ and
$G=(R^{\widetilde{n}},\widetilde{V}_{1},\widetilde{\delta}_{\lambda},\widetilde{\rho})$
be two Carnot groups where $\rho$ and $\widetilde{\rho}$
homogeneous norms endowed to $G$, $\widetilde{G}$ respectively.
Let $\Omega$ be a bounded open set of $G$. Let $1\leq
\alpha<\infty$. If $u\in{R^{1,\alpha}(\Omega,\widetilde{G})}$,
then

\begin{enumerate}
 \item $u$ is approximate Pansu differentiable a.e. in
  $\Omega$.
  Let $Du(p)$ be the approximate Pansu differential at
  $p\in{\Omega}$. The linear map $D_{h}u$ from $V_{1}$ to $\widetilde{V}_{1}$ can be extended to a
  homomorphism $\mathcal{D}_{u}(p)$
   of Lie algebras such that
 $$
  Du(p)=\widetilde{{\exp}}\circ\mathcal{D}_{u}(p)\circ\exp^{-1}.
 $$
   \item\label{cinfty} If $\widetilde{\rho}$ is a homogeneous
norm of the class of $C^{\infty}$ on $G\backslash\{0\}$, then for
a.e. $p\in{\Omega}$, $Du(p)$ is the Pansu differential of
  $u$ in the topology of $L^{\alpha}(\Omega)$. That is,
\begin{equation}\label{eq:Latopologyconvergence}
 \lim_{\epsilon\rightarrow 0}\int_{\rho(\omega)\leq
 1} \left(\widetilde{\rho}
((Du(p)(\omega))^{-1}\widetilde{\delta}_{\frac{1}{
 \epsilon}}(u(p)^{-1}u(p\delta_{\epsilon}\omega)))\right)^{\alpha}d\omega=0.
\end{equation}
\end{enumerate}
\end{theorem}

\subsection{Precompactness and the trace theorem for Sobolev
mappings}\label{sec:2-3}
 In this section we first give a
compactness theorem and then develop a trace theorem, which will
be needed in Section \ref{sec:4}. The trace theorem for Sobolev
mappings between metric spaces is delicate. In \cite{KS}, a
satisfactory trace theorem was developed for mappings in
$KS^{1,\alpha}(\Omega,M)$ when $\Omega$ is a Lipschitz Riemannian
domain and $M$ is a complete metric space. In the case $\Omega$ is
a sub-Riemannian domain, whether an analogue can be developed is
the problem we are going to investigate. Note that even for scalar
valued Sobolev functions the trace theorem is not trivial when the
domain is sub-Riemannian, see \cite{GN1} and \cite{DGN} for
extensive discussions. The difficulty to this problem is partly
due to the presence of characteristic points in the boundary of
domain. In this paper, we will not deal with the case when the
boundary of the domain possesses characteristic points. The
characteristic case will be investigated in a forthcoming paper.
\vskip10pt

We first have the following precompactness theorem. Since its
proof is standard (see \cite{Am}, Theorem 2.4 and \cite{KS},
Theorem 1.13), we omit it.
\begin{theorem}\label{thm:campactness}
Let $(R^{n},\Delta,d_{c})$ be a C-C space and $\Omega\subset
R^{n}$ be a bounded open set. Let $M$ be a separable complete
metric space with a (quasi-)metric $d$. Let $\alpha>1$. Assume
$\{u_{\mu}\}_{\mu=1}^{\infty}$ be a sequence of mappings in
$R^{1,\alpha}(\Omega,M)$ such that
$$
  \sup_{\mu}\left\{\int_{\Omega}d^{\alpha}(u_{\mu}(p),m_{0})dp+\int_{\Omega}g_{\mu}^{\alpha}(p)dp
\right\}\leq C
$$
where $0\leq g_{\mu}\in{L^{\alpha}(\Omega)}$ is a dominant
function of the horizontal derivatives of
$\theta_{m}^{\mu}(p):=d(m,u_{\mu}(p))$ for any $m\in{M}$, that is,
$|X\theta_{m}^{\mu}(p)|\leq g(p)$ a.e. $p\in{\Omega}$ for any
$m\in{M}$ (see Definition \ref{sec2:definition2}); $C>0$ is an
absolute constant; $m_{0}$ is a fixed point in $M$. Then there
exists a subsequence
$\{u_{\mu^{\prime}}\}_{\mu^{\prime}=1}^{\infty}$ of
$\{u_{\mu}\}_{\mu=1}^{\infty}$ and a mapping
$u\in{R^{1,\alpha}(\Omega,M)}$ such that
\begin{enumerate}
    \item
    $\lim_{\mu^{\prime}\rightarrow\infty}\int_{\Omega}d^{\alpha}(u_{\mu^{\prime}}(p),u(p))dp=0$;
    \item there exists a dominant function $0\leq g\in{L^{\alpha}(\Omega)}$
    of the horizontal derivatives of
     $\theta_{m}(p):=d(m,u(p))$ ($m\in{M}$) satisfies $\int_{\Omega}g^{\alpha}(p)dp
     \leq \lim\limits_{\mu^{\prime}\rightarrow0}\int_{\Omega}g_{\mu^{\prime}}^{\alpha}(p)dp\leq
     C$.
\end{enumerate}
\end{theorem}
Let $(R^{n},\Delta,d_{c})$ be a C-C space and $\Omega\subset
R^{n}$ be a $C^{2}$ smooth bounded domain whose boundary does not
possess characteristic points. We recall that \textit{a
characteristic point} $p\in{\partial\Omega}$ is a point where the
tangent space $T_{p}\partial\Omega$ contains the horizontal space
$\Delta_{p}$. Let $\vec{n}$ be the unit Euclidean exterior normal
vector field of $\partial\Omega$. Since $\Omega$ is $C^{2}$, there
exists a neighborhood $\widetilde{U}$ of $\partial\Omega$ such
that the signed distance function
$$d_{e}(p)=
  \begin{cases}
  -\textrm{dist}(p,\partial\Omega):=\inf\limits_{q\in{\partial\Omega}}|p-q| & \text{if } p\in{U\cap\Omega}, \\
    \textrm{dist}(p,\partial\Omega) & \text{if } p\in{U\cap\Omega^{c}}.
  \end{cases}
$$
is a defining function  of $\Omega$ (near the boundary), that is
$d_{e}$ is $C^{2}$ in $\widetilde{U}$ and $\vec{n}=\nabla d_{e}$
where $\vec{n}$ is the unit Euclidean exterior normal vector
fields in $\partial\Omega$. Since we have assumed that
$\partial\Omega$ is not characteristic, there exists a constant
$0<\widetilde{\rho}\leq 1$ such that  the horizontal transverse
vector field $Z(p)=Xd_{e}(p)=\sum_{i=1}^{k}<X_{i},\vec{n}>X_{i}$
satisfies
$$
|Z(p)|=\|Z(p)\|_{<\cdot,\cdot>_{c}}\geq
\widetilde{\rho}\quad\textrm{for any } p\in{U}\subset\widetilde{U}
$$
where $U$ is a neighborhood of $\partial\Omega$. The horizontal
transverse vector field $Z$ induces a fibration $\Gamma_{Z}$ of
$U\cap \Omega$, that is,
$\Gamma_{Z}=\{\gamma(t)=\exp_{p}(tZ):[0,t_{0}]\rightarrow \Omega,
p\in{\partial\Omega}\}$. Then $\gamma_{p}(t)$ satisfies
\begin{equation}\label{lipschitzdomain}
\begin{array}{ccl}
   &\gamma_{p}(t)\in{\Omega}& \textrm{ if } \quad0<t<t_{0},  \\
  &\gamma_{p}(t)\notin{\Omega}&\textrm{ if }\quad-t_{0}<t<0,\\
  &|d_{e}(\gamma_{p}(t_{1}))-d_{e}(\gamma_{p}(t_{2}))|>\rho|t_{1}-t_{2}|&\textrm{ if
  }\quad|t_{1}|, |t_{2}|<t_{0}
\end{array}
\end{equation}
for some choice of constants $\rho>0$ and $t_{0}>0$.

We recall the definition of the measure $d\gamma$ on $\Gamma_{Z}$,
$d\gamma=\mathcal{F}_{-t}i(Z)dv$, where $\mathcal{F}_{t}$ is the
Jacobian of the flow $\exp_{p}(tZ),p\in{\Omega}$ and $dv$ is the
standard volume form of $\Omega$. Since $Z$ is transversal to
$\partial\Omega$, the area form $d\sigma$ of $\partial\Omega$, up
to a normalization, is $i(Z)dv$ where $Z$ is understood as the
restriction on $\partial\Omega$ of $Z$. Note that
$\mathcal{F}_{-t}$ is always bounded in $U$.

 Let $M$ be a separable metric space
with a metric $d$. Let
$u\in{R^{1,\alpha}(\Omega,M)}(\alpha\geq1)$. We define the trace
$Tu\in{L^{\alpha}(\partial\Omega,M)}$ of $u$ on $\partial\Omega$
as follows.  By Theorem \ref{sec2:description} there exists a
representative $\widetilde{u}$ of $u$ such that $\widetilde{u}$ is
absolutely continuous on $d\gamma$ almost all curves in
$\Gamma_{Z}$, that is, $\widetilde{u}$ is absolutely continuous on
$\gamma_{p}(t)=\exp_{p}(tZ)$ ($0<t\leq t_{0}$) for $d\sigma$
almost all $p\in{\partial\Omega}$. Thus the map
$$
Tu(p)=\lim_{t\rightarrow 0^{+}}\widetilde{u}(\gamma_{p}(t))
$$
can be defined for a.e. $p\in{\partial\Omega}$. Furthermore from
the proof of Theorem \ref{sec2:description} (see (4.1) in Page 641
of \cite{Vo1}) and using H\"{o}lder inequality we have
\begin{equation}\label{trace1}
 d^{\alpha}(Tu(p),\widetilde{u}(\gamma_{p}(t)))
 \leq t^{\alpha-1}\int_{[p,\gamma_{p}(t)]}g^{\alpha}ds
\end{equation}
for a.e. $p\in{\partial\Omega}$ where $0\leq
g\in{L^{\alpha}(\Omega)}$. Integrating \eqref{trace1} with respect
to $p$ we infer
\begin{align}
 \int_{\partial\Omega}d^{\alpha}(Tu(p),\widetilde{u}(\gamma_{p}(t)))d\sigma(p)
 &\leq t^{\alpha-1}\int_{\partial\Omega}\int_{[p,\gamma_{p}(t)]}g^{\alpha}dsd\sigma(p)\notag\\
 &\leq C
 t^{\alpha-1}\int_{\Omega_{t}^{C}}g^{\alpha}dv\label{trace2}
\end{align}
where $C$ is a constant independent of $t$ and
$\Omega_{\epsilon}^{C}$ denotes the set of points in $\Omega$
whose C-C distance to $\partial\Omega$ is at most $\epsilon$.
Since $\widetilde{u}\in{L^{\alpha}(\Omega,M)}$, by the Fubini's
theorem the maps $\widetilde{u}(\gamma_{p}(t))$ are in
$L^{\alpha}(\partial\Omega,M)$ for almost all $t\in[0,t_{0}]$. We
conclude from \eqref{trace2} that the trace map $T(u)$ is the
$L^{\alpha}(\partial\Omega,M)$ limit of the maps
$\widetilde{u}(\gamma_{p}(t))$ as $t\rightarrow0$, so is itself an
$L^{\alpha}$ map. Since $T(u)$ is the $L^{\alpha}$ limit of almost
all of the maps $\widetilde{u}(\gamma_{p}(t))$, as $t\rightarrow
0$, $T(u)$ is independent of the choice of the representative of
$u$.

Thus we have proven
\begin{proposition}\label{pro:trace}
Let $(R^{n},\Delta,d_{c})$ be a C-C space and $\Omega\subset
R^{n}$ be a $C^{2}$ smooth bounded domain whose boundary
$\partial\Omega$ does not possess characteristic points. Let $M$
be a separable metric space with a (quasi-)metric $d$. Assume
$u\in{R^{1,\alpha}(\Omega,M)}$, $\alpha\geq 1$. Then the trace
$T(u)$ is well defined and $Tu\in{L^{\alpha}(\partial\Omega,M)}$.
\end{proposition}
\begin{remark}
In Proposition \ref{pro:trace} the noncharacteristic condition is
restrictive in the sense that ``most" smooth bounded domains in a
C-C space are characteristic. For examples of noncharacteristic
smooth bounded domains we refer the reader to \cite{DGN}.
\end{remark}
The following lemma can be easily deduced from Theorem
\ref{sec1:equivalence}, Theorem \ref{sec2:equivalence}, Corollary
\ref{energyequivalence} and Corollary 1.6.3 in \cite{KS}.
\begin{lemma}
Let $(R^{n},\Delta,d_{c})$ be a C-C space and $\Omega\subset
R^{n}$ be a bounded open set. Let $M$ be a separable metric space
with a (quasi-)metric $d$ and let $\alpha>1$. Assume
$u,v\in{R^{1,\alpha}(\Omega,M)}$ with dominant functions  $g_{u}$,
$g_{v}$ respectively, then $d(u,v)\in{W_{X}^{1,\alpha}(\Omega)}$
and
\begin{equation}\label{eq:difference}
   \|Xd(u,v)\|_{L^{\alpha}(\Omega)}\leq C(\|g_{u}\|_{L^{\alpha}(\Omega)}+\|g_{v}\|_{L^{\alpha}(\Omega)})
\end{equation}
for some constant $C>0$.
\end{lemma}
\begin{theorem}\label{tracecharacterization}
Let $(R^{n},\Delta,d_{c})$ be a C-C space and $\Omega\subset
R^{n}$ be a $C^{2}$ smooth bounded domain whose boundary
$\partial\Omega$ does not possess characteristic points. Let $M$
be a separable complete metric space with a (quasi-)metric $d$ and
let $\alpha>1$. If the sequence
$\{u_{\mu}\}_{\mu=1}^{\infty}\subset R^{1,\alpha}(\Omega,M)$ has a
sequence of dominant functions
$\{g_{\mu}\}_{\mu=1}^{\infty}\subset L^{\alpha}(\Omega)$ (that is,
$g_{\mu}$ is a dominant function for $u_{\mu}$) such that
$\{\|g_{\mu}\|^{\alpha}_{L^{\alpha}(\Omega)}\}_{\mu=1}^{\infty}$
has uniform bound, and if $\{u_{\mu}\}$ converges in
$L^{\alpha}(\Omega,M)$ to a mapping $u:\Omega\rightarrow M$, then
the trace functions of $u_{\mu}$ converge in
$L^{\alpha}(\partial\Omega,M)$ to the trace of $u$. Two mapping
$u,v\in{R^{1,\alpha}(\Omega,M)}$ have the same trace if and only
if $d(u,v)\in{W_{X}^{1,\alpha}(\Omega)}$ has trace zero.
\end{theorem}
\begin{proof} It follows from Theorem \ref{thm:campactness} that the
$L^{\alpha}(\Omega,M)$ limit map $u$ belongs to
$R^{1,\alpha}(\Omega,M)$. Since \eqref{trace2} holds, Theorem
\ref{tracecharacterization} follows almost verbatim from the
arguments in the proof of Theorem 1.12.2 in \cite{KS}( for the
existence of $d_{c}$-Lipschitz cut-off functions see \cite{GN1}).

\end{proof}
\section{Energy of Korevaar-Schoen}\label{sec:3}
This section and Section \ref{sec:4} are devoted to making a
choice of a reasonable energy which should be natural and
compatible to the structures of the considered C-C spaces. Since
the energy of Korevaar-Schoen has been extensively studied, a
natural question is that whether it is the one we expected. When
the target does not possess any curvature bound in the sense of
Alexandrov, we want the energy to be of ``good" form, for example,
it is a Dirichlet integral. Unfortunately the energy of
Korevaar-Schoen is not of the form of the Dirichlet integral,
though it can be represented by an integral (see
\eqref{korevaar-schoenenergy}). By now we can not prove or
disprove that the energy of Korevaar-Schoen is lower
semicontinuous with respect to some topology. Note that C-C spaces
may not possess ``measure contraction property" which Riemannian
manifolds possess (see \cite{Sturm}). Thus we can not adopt the
idea in \cite{KS} and \cite{Sturm}.

\begin{theorem}\label{thm:energychara}
Let $G=(R^{n},V_{1},\delta_{\lambda},\rho)$ and
$G=(R^{\widetilde{n}},\widetilde{V}_{1},\widetilde{\delta}_{\lambda},\widetilde{\rho})$
be two Carnot groups where $\rho$ and $\widetilde{\rho}$ are
homogeneous norms endowed to $G$ and $\widetilde{G}$ respectively.
Let $\Omega$ be a bounded open set of $G$. Let
$\alpha\in{(1,\infty)}$. If $\widetilde{\rho}$ is of the class of
$C^{\infty}$ on $\widetilde{G}\backslash\{0\}$ and
$u\in{KS^{1,\alpha}(\Omega,\widetilde{G})}$, then the energy of
Korevaar-Schoen can be written as:
\begin{equation}\label{korevaar-schoenenergy}
  E^{\alpha}(u,\Omega)=\int_{\Omega}\dashint_{B_{\rho}(0,1)}
  \left(\widetilde{\rho}(Du(p)(\omega))\right)^{\alpha}d\omega dp.
\end{equation}
where $B_{\rho}(0,1)=\{\omega:\rho(\omega)\leq 1\}$ and
$Du(p):G\rightarrow \widetilde{G}$ is the approximate Pansu
derivative of $u$ at $p$.
\end{theorem}
\begin{proof} We abuse the notation
$\rho(p,q):=d_{\rho}(p,q)=\rho(p^{-1}q)$. By Theorem
\ref{pansuofSobolev} $u$ is approximately Pansu differentiable
a.e. $p\in{\Omega}$. Fix $p\in{\Omega}$ at which $u$ is
approximately Pansu differentiable in $\Omega$. Recalling that the
Lebesgue measure $\mathcal{L}^{n}$ in $R^{n}$ is the Haar measure
of $G$ and the definition of homogeneous norms, by a change of
variables we have
\begin{align}\label{ksenergy}
  e_{\epsilon}(p;u) &= \dashint_{B_{\rho}(p,\epsilon)}\left(\frac{
\widetilde{\rho}(u(p),u(q))}{\epsilon}\right)^{\alpha}dq\notag\\
   & =\dashint_{B_{\rho}(0,1)}\left(\widetilde{
   \rho}(\widetilde{\delta}_{\frac{1}{\epsilon}}(u(p)^{-1}u(p\delta_{\epsilon}
   \omega)))\right)^{\alpha}d\omega
\end{align}

By Theorem \ref{sec2:equivalence}, we have
$u\in{R^{1,\alpha}(\Omega,\widetilde{G})}$. Now we can use
\eqref{eq:Latopologyconvergence} to deduce
\begin{equation}\label{density}
\lim_{\epsilon\rightarrow0}e_{\epsilon}(p;u)=\dashint_{B_{\rho}(0,1)}
  \left(\widetilde{\rho}(Du(p)\omega)\right)^{\alpha}d\omega.
\end{equation}
In fact, by a well known inequality
\begin{equation}\label{elementary1}
||a+b|^{\alpha}-|b|^{\alpha}|\leq
C(\delta)|a|^{\alpha}+\delta|b|^{\alpha}
\end{equation}
where $a,b\in{R}$, $\delta>0$ and $C(\delta)$ only depends on
$\delta$ and $\alpha$, we obtain
\begin{align}
 &\left(\widetilde{\rho}(\widetilde{\delta}_{\frac{1}{\epsilon}}(u(p)^{-1}u(p\delta_{\epsilon}\omega)))\right)^{\alpha}
 -\left(\widetilde{\rho}(Du(p)\omega)\right)^{\alpha}\notag\\
  &=\left(\widetilde{\rho}(\widetilde{\delta}_{\frac{1}{\epsilon}}(u(p)^{-1}u(p\delta_{\epsilon}\omega)))-
  \widetilde{\rho}(Du(p)\omega)+\widetilde{\rho}(Du(p)\omega)\right)^{\alpha}
 -\left(\widetilde{\rho}(Du(p)\omega)\right)^{\alpha}\notag\\
 &\leq\delta\left(\widetilde{\rho}(Du(p)\omega)\right)^{\alpha} +C(\delta)\left(\widetilde{\rho}(\widetilde{\delta}_
 {\frac{1}{\epsilon}}(u(p)^{-1}u(p\delta_{\epsilon}\omega)))-
  \widetilde{\rho}(Du(p)\omega)\right)^{\alpha}\label{elementary2}\\
  &\leq\delta\left(\widetilde{\rho}(Du(p)\omega)\right)^{\alpha}+C(\delta)C^{\prime}
  \left(\widetilde{\rho}((Du(p)\omega)^{-1}\widetilde{\delta}_
 {\frac{1}{\epsilon}}(u(p)^{-1}u(p\delta_{\epsilon}\omega)))
  \right)^{\alpha}\label{triangleine}
\end{align}
where \eqref{elementary2} is from \eqref{elementary1} and in
\eqref{triangleine} we have used the quasi-triangle inequality
property of $\widetilde{\rho}$. Thus \eqref{density} follows from
\eqref{ksenergy}, \eqref{triangleine},
\eqref{eq:Latopologyconvergence} and the arbitrariness of $\delta$
in \eqref{triangleine}.

On the other hand, since by Theorem \ref{sec2:equivalence}
$u\in{H^{1,\alpha}(\Omega,\widetilde{G})}$, there exists $0\leq
g\in{L^{\alpha}(\Omega)}$ such that
\begin{equation}\label{uinhajlasz}
  \widetilde{\rho}(u(p),u(q))\leq \rho(p,q)(g(p)+g(q))
\end{equation}
holds for a.e. $p,q\in{\Omega}$. Let
$\varphi\in{C_{c}(\Omega,[0,1])}$ and
$\epsilon<\textrm{dist}_{\rho}(\textrm{supp}\varphi,\partial\Omega)$.
Assume
\begin{equation}\label{control3}
 F_{\epsilon}(u,p):=\varphi(p)e_{\epsilon}(p;u),
\end{equation}
 then by \eqref{uinhajlasz} we have
\begin{equation}\label{compare}
 F_{\epsilon}(u,p)\leq G_{\epsilon}(u,p):=C\varphi(p)
|g(p)|^{\alpha}+C\varphi(p)\dashint_{B_{\rho}(p,\epsilon)}|g(q)|^{\alpha}dq
\end{equation}
 for a.e. $p\in{\Omega}$ where $C$ only depends on $\alpha$.
 Since $g\in{L^{\alpha}(\Omega)}$, by Lebesgue differentiation theorem
 (see e.g. \cite{Heinonen}, Chapter 2)
\begin{equation}\label{compare1}
 \lim_{\epsilon\rightarrow0}G_{\epsilon}(u,p)=G(u,p):=2C\varphi(p)|g(p)|^{\alpha}\quad
 \textrm{a.e. } p\in{\Omega}.
 \end{equation}
From
\begin{align}
  \int_{\Omega}G_{\epsilon}(u,p)dp
  &=C\int_{\Omega}\varphi(p)|g(p)|^{\alpha}dp+
  C\int_{\Omega}\varphi(p)\dashint_{B_{\rho}(p,\epsilon)}|g(q)|^{\alpha}dqdp\notag\\
  &
  =C\int_{\Omega}\varphi(p)|g(p)|^{\alpha}dp+C
  \int_{\Omega_{\epsilon}}\varphi(p)\dashint_{B_{\rho(0,1)}}|g(p\delta_{\epsilon}q^{\prime})|^{\alpha}dq^{\prime}dp\label{6}\\
  &=C\int_{\Omega}\varphi(p)|g(p)|^{\alpha}dp+C\dashint_{B_{\rho(0,1)}}\int_{\Omega_{\epsilon}}
  \varphi(p)|g(p\delta_{\epsilon}q^{\prime})|^{\alpha}dpdq^{\prime}\label{7}\\
  &=C\int_{\Omega}\varphi(p)|g(p)|^{\alpha}dp+C\dashint_{B_{\rho(0,1)}}\int_{\Omega}
  \varphi(p^{\prime}\delta_{\epsilon}q^{\prime^{-1}})|g(p^{\prime})|^{\alpha}dp^{\prime}dq^{\prime}\label{8}
 \end{align}
where in \eqref{6} we have made the change of variables
$q^{\prime}=\delta_{\frac{1}{\epsilon}}(p^{-1}q)$, in \eqref{7}
used the Fubini's Theorem, in \eqref{8} used the change of
variables $p^{\prime}=p\delta_{\epsilon}q^{\prime}$,
$$
\varphi(p^{\prime}\delta_{\epsilon}q^{\prime^{-1}})|g(p^{\prime})|^{\alpha}\leq
|g(p^{\prime})|^{\alpha}
$$
and
$$
 \int_{\Omega}
  \varphi(p^{\prime}\delta_{\epsilon}q^{\prime^{-1}})|g(p^{\prime})|^{\alpha}dp^{\prime}dq^{\prime}\leq
  \|g\|^{\alpha}_{L^{\alpha}(\Omega)},
$$
by dominated convergence theorem we infer that
\begin{align}
 \lim_{\epsilon\rightarrow0}\int_{\Omega}G_{\epsilon}(u,p)dp
  &= C\int_{\Omega}\varphi(p)|g(p)|^{\alpha}dp+C\dashint_{B_{\rho(0,1)}}
  \int_{\Omega}
  \lim_{\epsilon\rightarrow0}\varphi(p^{\prime}\delta_{\epsilon}q^{\prime^{-1}})
  |g(p^{\prime})|^{\alpha}dp^{\prime}dq^{\prime}\notag\\
  &=2C\int_{\Omega}\varphi(p)|g(p)|^{\alpha}dp\notag\\
  &=\int_{\Omega}G(u,p)dp\label{generalized}
\end{align}
since $\varphi\in{C_{c}(\Omega,[0,1])}$. From \eqref{density},
\eqref{control3}, \eqref{compare}, \eqref{compare1},
\eqref{generalized} and a variant dominated convergence theorem
(see e.g. \cite{EG}, p21), we have
$$
\lim_{\epsilon\rightarrow0}\int_{\Omega}F_{\epsilon}(u,p)dp=\int_{\Omega}\varphi(p)\dashint_{B_{\rho}(0,1)}
  \left(\widetilde{\rho}(Du(p)\omega)\right)^{\alpha}d\omega dp.
$$
for any $\varphi\in{C_{c}(\Omega,[0,1])}$. Consequently
\eqref{ksenergy} follows.
\end{proof}

In Theorem \ref{ksenergy} we give a representation of the energy
of Korevaar-Schoen for mappings in
$KS^{1,\alpha}(\Omega,\widetilde{G})$ from a Carnot group to
another Carnot group. One may ask whether $E^{\alpha}(u,\Omega)$
is lower semicontinuous with respect to some topology of $
KS^{1,\alpha}(\Omega,\widetilde{G})$. When $\Omega$ is a smooth
Riemannian domain and $\widetilde{G}$ is a metric space with a
metric $d$, in \cite{KS} Korevaar-Schoen proved
$E^{\alpha}(u,\Omega)$ is lower semicontinuous with respect to the
topology of $L^{\alpha}(\Omega,\widetilde{G})$, that is, if
$\sup_{\mu} E^{\alpha}(u_{\mu},\Omega)<\infty$ and
$\lim_{\mu\rightarrow\infty}\int_{\Omega}d(u_{\mu}(p),u(p))^{\alpha}dp=0
$, then $ E^{\alpha}(u,\Omega)\leq
\liminf_{\mu\rightarrow\infty}E^{\alpha}(u_{\mu},\Omega)$. This
lower semicontinuity property is based on a subpartitional lemma
for the approximate energies $E_{\epsilon}^{\alpha}(\varphi;u)$
(see \eqref{approximateenergy}) which plays a fundamental role in
the whole story of \cite{KS}. The idea in \cite{KS} of
constructing  Sobolev mappings between metric spaces has been used
and generalized by several authors, see \cite{Sturm}, \cite{KS1},
\cite{KS2} and \cite{EF}. Whether or not a subpartitional lemma
for the approximate energies holds depends on the metric property
of the domain space and is independent of the target. Sturm in
\cite{Sturm} proposed a type of metric spaces which possess so
called a measure contraction property(MCP). The class of MCP
spaces includes Lipschitz Riemannian spaces. Sturm proved that a
subpartitional lemma for the approximate energies holds in MCP
spaces. Now a natural question, which has independent interests,
is whether or not C-C spaces are MCP spaces (or SMCPBG, GMCP in
the sense of \cite{KS1}, \cite{KS2}). In the following, for
Heisenberg group we will prove that the Jacobian of the change of
variables along C-C geodesics is not what we expected. \vskip10pt

Let's first recall some fundamental facts about C-C geodesics in
Heisenberg groups $H^{m}$(see  \cite{AR} or \cite{TY}).
\begin{lemma}\label{ccgeodesic}
Let $g_{0}=(x_{0},y_{0},t_{0})\neq 0$ be a point in $H^{m}$. We
have
\begin{enumerate}
  \item if $x_{0}^{2}+y_{0}^{2}\neq 0$, then there exists a unique C-C geodesic connecting $0$ to $g_{0}$.
  \item otherwise, there exist infinitely many  C-C
  geodesics connecting $0$ to $g_{0}$.
\end{enumerate}
Moreover, let $\gamma(s)=(x(s),y(s),t(s))(0\leq s\leq 1)$ be any
C-C geodesic connecting $0$ to $g_{0}$, we have
$$
\left\{
\begin{array}{rl}
x_{i}(s)&=\frac{A_{i}(\cos(s\phi\rho) -1)+B_{i}\sin(s\phi\rho)}{\phi},\quad i=1,\cdots,m,\\
y_{i}(s)&=\frac{B_{i}(\cos(s\phi\rho) -1)-A_{i}\sin(s\phi\rho)}{\phi},\quad i=1,\cdots,m,\\
t(s)&=2\frac{s\phi\rho-\sin (s\phi\rho)}{\phi^{2}},
\end{array}
\right.
$$
where $\tau=\phi\rho\in{[-2\pi,2\pi]}$ is  the unique solution in
$[-2\pi,2\pi]$ of the equation
\begin{equation}\label{tau-Heisenberg}
\frac{1-\cos\tau}{\tau-\sin\tau}=\frac{|x_{0}|^{2}+|y_{0}|^{2}}{t_{0}}
\end{equation}
with
$$
  \begin{cases}
   \tau=0& \text{if }t_{0}=0,\\
    |\tau|=2\pi&\text{if }|x_{0}|^{2}+|y_{0}|^{2}=0,\\
  \tau\in{(0,2\pi)}&\text{if }t_{0}>0,\\
  \tau\in{(-2\pi,0)}&\text{otherwise};
  \end{cases}
  $$
 $\rho=d_{c}(0,g_{0})$ is the arc length of $\gamma$ determined by
 \begin{eqnarray*}
\rho=&\sqrt{\frac{\tau^{2}t_{0}}{2(\tau-\sin\tau)}},\quad&\textrm{
if }t_{0}\neq 0,
\\
\rho=&\sqrt{|x_{0}|^{2}+|y_{0}|^{2}},\quad&\textrm{ if }t_{0}=0;
\end{eqnarray*}
if $|x_{0}|^{2}+|y_{0}|^{2}\neq0$
$\{A_{1},\cdots,A_{m},B_{1},\cdots,B_{m}\}$ is subject to
$$
\begin{cases}
   \sum_{i=1}^{n}(A_{i}^{2}+B_{i}^{2})=1, \\
   {x_{0}}_{i}=\frac{A_{i}(\cos(\phi\rho) -1)+B_{i}\sin(\phi\rho)}{\phi},\quad i=1,\cdots,m,\\
   {y_{0}}_{i}=\frac{B_{i}(\cos(\phi\rho) -1)-A_{i}\sin(\phi\rho)}{\phi},\quad
   i=1,\cdots,m;
  \end{cases}
$$ if $|x_{0}|^{2}+|y_{0}|^{2}=0$ then
$\{A_{1},\cdots,A_{m},B_{1},\cdots,B_{m}\}$ is only subject to
$$
\sum_{i=1}^{m}(A_{i}^{2}+B_{i}^{2})=1.
$$
\end{lemma}
\begin{remark}\label{lefttranslationgeodesic}
By the left-invariant property of the C-C metric we easily deduce
that $\gamma_{p_{0},p}(s)$ is a C-C geodesic connecting $p_{0}$ to
$p$ if and only if $p_{0}^{-1}\gamma_{p_{0},p}(s)$ is a C-C
geodesic connecting $0$ to $p_{0}^{-1}p$, that is,
$\gamma_{p_{0},p}(s)=p_{0}\gamma_{p_{0}^{-1}p}(s)$ where
$\gamma_{p_{0}^{-1}p}(s)$ is the C-C geodesic connecting $0$ to
$p_{0}^{-1}p$.
\end{remark}
To simplify some computation we fix $m=1$. Let
$H_{*}^{1}=H^{1}\backslash \{(0,0,t):t\in{R}\}$. Set
$S=\{(\theta,\phi,\rho):0\leq \theta<2\pi,|\phi\rho|\leq
2\pi,\rho\geq0\}$, $S_{*}=\{(\theta,\phi,\rho):0\leq
\theta<2\pi,|\phi\rho|<2\pi,\rho\geq0\}$ and
$\mathcal{A}(\theta,\phi,\rho):S\rightarrow H^{1}$ by
$\mathcal{A}(\theta,\phi,\rho)=(x(\theta,\phi,\rho),y(\theta,\phi,\rho),t(\theta,\phi,\rho))$,
where
\begin{equation}\label{parametrization}
\left\{\begin{array}{rl}
  x(\theta,\phi,\rho) & = \frac{\cos\theta(\cos(\phi\rho) -1)+\sin\theta\sin(\phi\rho)}{\phi}\\
  y(\theta,\phi,\rho) & =\frac{\sin\theta(\cos(\phi\rho)
  -1)-\cos\theta\sin(\phi\rho)}{\phi}\\
  t(\theta,\phi,\rho) &=2\frac{\phi\rho-\sin (\phi\rho)}{\phi^{2}}
\end{array}\right.
\end{equation}
By Lemma \ref{ccgeodesic}, we know that the map
$\mathcal{A}:S_{*}\rightarrow H^{1}_{*}$ is bijective and equation
\eqref{parametrization} parameterizes $\partial
B_{c}(0,\rho)=\{p\in{H^{1}}:d_{c}(0,p)=\rho\}$. We can compute the
Jacobian of $\mathcal{A}$ by
\begin{equation}\label{jacobian1}
  \det J\mathcal{A}(\theta,\phi,\rho)=\det\left(\frac{\partial(x,y,t)}{\partial(\theta,\phi,\rho)}\right)=
  4\frac{\phi\rho\sin(\phi\rho)-2(1-\cos\phi\rho)}{\phi^{4}}
\end{equation}
Let $\bar{s}\in[0,1]$ and $p=(x,y,t)\in{H^{1}_{*}}$, we consider
the Jacobian of the map of changing variables
$\mathcal{B}_{0}^{\bar{s}}(p):p=(x,y,t)\rightarrow
p^{\prime}=(x^{\prime},y^{\prime},t^{\prime})=\gamma_{0,p}(\bar{s})$
where $\gamma_{0,p}(s)(0\leq s\leq 1)$ is the C-C geodesic joining
0 and $p$. Since $p\in{H^{1}_{*}}$, we can parameterize $p$ by
$(\theta,\phi,\rho)$ through equation \eqref{parametrization} and
$(x^{\prime},y^{\prime},t^{\prime})=\mathcal{A}(\theta,\phi,\bar{s}\rho)$.
Now we can compute the Jacobian of $\mathcal{B}_{0}^{\bar{s}}$ by
\begin{align}
  \det J\mathcal{B}_{0}^{\bar{s}}(x,y,t)& =\det\left(\frac{\partial(x^{\prime},y^{\prime},t^{\prime})}
  {\partial(x,y,t)}\right)\notag\\
   &=
   \det\left(\frac{\partial(x^{\prime},y^{\prime},t^{\prime})}{\partial(\theta,\phi,\bar{s}\rho)}\right)
   \det\left(\frac{\partial(\theta,\phi,\bar{s}\rho)}{\partial(\theta,\phi,\rho)}\right)
   \det\left(\frac{\partial(\theta,\phi,\rho)}{\partial(x,y,t)}\right)\notag\\
   &=\det J\mathcal{A}(\theta,\phi,\bar{s}\rho)\,\bar{s}\,\left(\det J\mathcal{A}(\theta,\phi,\rho)\right)^{-1}\notag\\
   &=\bar{s}\frac{\bar{s}\phi\rho\sin(\bar{s}\phi\rho)-2(1-\cos(\bar{s}\phi\rho))}
   {\phi\rho\sin(\phi\rho)-2(1-\cos(\phi\rho))}\label{jacobian2}
\end{align}
where we have used \eqref{jacobian1}. From \eqref{tau-Heisenberg}
we have if $t_{0}\rightarrow0$, then $\tau=\phi\rho\rightarrow0$.
Thus from \eqref{jacobian2} we obtain
\begin{equation}\label{jocabian3}
  \lim_{t\rightarrow0}\det
  J\mathcal{B}_{0}^{\bar{s}}(x,y,t)=\lim_{\tau\rightarrow0}\bar{s}\frac{\bar{s}
  \tau\sin(\bar{s}\tau)-2(1-\cos(\bar{s}\tau))}
   {\tau\sin\tau-2(1-\cos\tau)}
   =\bar{s}^{5}
\end{equation}

In general case for any
$p_{0}=(x_{0},y_{0},t_{0})\in{H^{1}},\bar{s}\in{[0,1]}$, we define
the map
$\mathcal{B}_{p_{0}}^{\bar{s}}(p)=\gamma_{p_{0},p}(\bar{s})=p_{0}\gamma_{p_{0}^{-1}p}(\bar{s})$
where
$p\in{{}_{*}H^{1}_{{p_{0}}}=\{g\in{H^{1}}:p_{0}^{-1}g\in{H^{1}_{*}}\}}$,
$\gamma_{p_{0},p}(s)$ is the C-C geodesic connecting $p_{0}$ to
$p$ and $\gamma_{p_{0}^{-1}p}(s)$ is the C-C geodesic connecting
$0$ to $p_{0}^{-1}p$ (see Remark \ref{lefttranslationgeodesic}).
Since the Jacobian of left translation is 1, from
\eqref{jacobian2} we can easily infer that the Jacobian of
$\mathcal{B}_{p_{0}}^{\bar{s}}$ is
\begin{equation}\label{jacobian6}
  \det J\mathcal{B}_{p_{0}}^{\bar{s}}(p) =\bar{s}\frac{\bar{s}\phi\rho\sin(\bar{s}\phi\rho)-2(1-\cos(\bar{s}\phi\rho))}
   {\phi\rho\sin(\phi\rho)-2(1-\cos(\phi\rho))}
\end{equation}
where $(\theta,\phi,\rho)$ parameterizes the point $p_{0}^{-1}p$
through \eqref{parametrization}. Thus we deduce from
\eqref{jacobian6} and \eqref{tau-Heisenberg} that
\begin{equation}\label{jacobian5}
  \lim_{t-t_{0}+2(x_{0}y-y_{0}x)\rightarrow0}\det J\mathcal{B}_{p_{0}}^{\bar{s}}(p)=\bar{s}^{5}
\end{equation}
where $p=(x,y,t)\in{{}_{*}H^{1}_{p_{0}}}$. \vskip10pt

It is well known that in a $m$ dimensional smooth  Riemannian
manifold $M$ with Riemannian metric $g$ the Jacobian of the map of
changing variables along Riemannian geodesics can be well
estimated. More precisely, if $\epsilon>0$ is sufficiently small
such that the geodesic ball $B_{d_{g}}(p_{0},\epsilon)$ is in  a
normal coordinate neighborhood of $p_{0}$, then $\det
J\mathcal{B}_{p_{0}}^{\bar{s}}(p)\geq
C\frac{\bar{s}^{m}}{1+o(\epsilon)}$ for any
$p\in{B_{g}(p_{0},\epsilon)}$ and $\bar{s}\in{[0,1]}$ where the
map $\mathcal{B}_{p_{0}}^{\bar{s}}(p)$ is defined similarly as
above and $C>0$ is a constant dependent on the Riemannian metric
(Ricci curvature). In the case of Heisenberg group, since the
Hausdorff dimension of $H^{m}$ is $Q=2m+2$, one would like to
guess that
\begin{equation}\label{jacobian4}
 \det J\mathcal{B}_{p_0}^{\bar{s}}(p)\geq C\bar{s}^{Q}
 \quad\textrm{in a neighborhood  of }p_{0}
\end{equation}
 for some constant $C$. We remark that
 if \eqref{jacobian4} was true, then several problems in analysis on Heisenberg
 groups could be solved by standard methods, for
 example to prove an inequality conjectured by \cite{AM} and
 to prove the semicontinuity of the energy of Korevaar-Schoen when
 the domain space is a Heisenberg group by repeating  the story of
 \cite{KS} or \cite{Sturm}.
 Unfortunately as we have shown above, \eqref{jacobian4} is impossible to hold.
\section{Horizontal energy and existence of minimizers}\label{sec:4}

As we have indicated in the Introduction, the concept of the
Horizontal energy is a natural generalization of the ordinary
energy for mappings between Riemannian manifolds.
\begin{definition}\label{def:horizontalenergy}
Let $(R^{n},\Delta,d_{c})$ and
$M=(R^{\widetilde{n}},\widetilde{\Delta},\widetilde{d}_{c})$ be
two C-C space spaces. Let $1\leq\alpha<\infty$ and $\Omega$ be a
bounded open set of $R^{n}$. Let $u\in{R^{1,\alpha}}(\Omega,M)$,
we call the following quantity
\begin{equation}
HE^{\alpha}(u,\Omega)=\int_{\Omega}\left(\sum_{i=1}^{k}\|X_{i}u(p)\|^{2}_
{\widetilde{<\cdot,\cdot>}_{c}}\right)^{\frac{\alpha}{2}}dp
\end{equation}
the $\alpha$-horizontal energy of $u$.
\end{definition}
Note that if $M$ is of Carnot type, from Remark \ref{remark:cha}
we have
$$
HE^{\alpha}(u,\Omega)=\int_{\Omega}\left(\sum_{j=1}^{k}
\sum_{i=l}^{\widetilde{k}}|X_{j}u^{i}(p)|^{2}\right)^{\frac{\alpha}{2}}dp.
$$

Let $(R^{n},\Delta,d_{c})$ be a C-C space and $\Omega\subset
R^{n}$ be a $C^{2}$ bounded open set whose boundary is
noncharacteristic with respect to $\Delta$. Let
$M=(R^{\widetilde{n}},\widetilde{\Delta},\widetilde{d}_{c})$ be
another C-C space and let $\alpha\geq 1$. Fix
$\phi\in{R^{1,\alpha}(\Omega,M)}$ and set
$$
R_{\phi}^{1,\alpha}(\Omega,M):=\{u\in{R^{1,\alpha}(\Omega,M):T(u)=T(\phi)}\},
$$
where $T(u)$ denotes the trace map of $u$, see Subsection
\ref{sec:2-3}.
 We consider the following Dirichlet problem of minimizing
$\alpha$-horizontal energy among all mappings in
$R^{1,\alpha}(\Omega,M)$ whose traces are equivalent to the trace
of $\phi$:
\begin{equation}\label{dirichletproblem}
   \textrm{find a }u\in{R_{\phi}^{1,\alpha}(\Omega,M)}\textrm{ such that
   }HE^{\alpha}(u,\Omega)=\inf_{v\in{R_{\phi}^{1,\alpha}(\Omega,M)}}HE^{\alpha}(v,\Omega).
\end{equation}
\begin{definition}\label{horizontalminimizer}
Any solution to Problem \eqref{dirichletproblem} is called a
horizontal energy minimizer.
\end{definition}
\begin{theorem}
Let $(R^{n},\Delta,d_{c})$ be a C-C space and $\Omega\subset
R^{n}$ be a $C^{2}$ bounded open set whose boundary is
noncharacteristic with respect to $\Delta$. Let
$M=(R^{\widetilde{n}},\widetilde{\Delta},d)$ is a C-C space of
Carnot type with a (quasi-)metric $d$ which is equivalent to
$\widetilde{d}_{c}$ and let $\alpha>1$. Then Problem
\eqref{dirichletproblem} has a solution.
\end{theorem}
\begin{proof} Let $\{u_{\mu}\}_{\mu=1}^{\infty}$ be a minimizing
sequence, that is,
\begin{equation}\label{minimizingsequence}
\lim_{\mu\rightarrow\infty}HE^{\alpha}(u_{\mu},\Omega)=
\inf_{v\in{R_{\phi}^{1,\alpha}(\Omega,M)}}HE^{\alpha}(v,\Omega):=C_{0}\leq
HE^{\alpha}(\phi,\Omega)<\infty.
\end{equation}
From Theorem \ref{sec2:characterization} and Remark
\ref{generalccspace} we easily get a sequence of dominant
functions $\{g_{\mu}\}_{\mu=1}^{\infty}$ of
$\{u_{\mu}\}_{\mu=1}^{\infty}$ such that
$\|g_{\mu}\|_{L^{\alpha}(\Omega)}\leq
C\,HE^{\alpha}(u_{\mu},\Omega)$ for any $\mu$ where $C$ is a
constant. From \eqref{minimizingsequence} the sequence
$\{g_{\mu}\}$ is uniformly bounded in $L^{\alpha}(\Omega)$. On the
other hand, since $T(u_{\mu})=T(\phi)$, by Theorem
\ref{tracecharacterization} we get
$d(u_{\mu},\phi)\in{W_{X}^{1,\alpha}(\Omega)}$ has trace zero for
any $\mu$. Applying the Poincar\'{e} inequality, (quasi-)triangle
inequality property of $d$ and \eqref{eq:difference} we get
\begin{align*}
  \int_{\Omega}d^{\alpha}(u_{\mu}(p),m_{0})dp &\leq
  C\left( \int_{\Omega}d^{\alpha}(u_{\mu}(p),\phi(p))dp+ \int_{\Omega}d^{\alpha}(m_{0},\phi)dp\right) \\
  \quad & \leq C_{1}\int_{\Omega}|Xd(u_{\mu},\phi)|^{\alpha}(p)dp+C\int_{\Omega}d^{\alpha}(m_{0},\phi)dp \\
  \quad & \leq C_{3}(\|g_{\mu}\|^{\alpha}_{L^{\alpha}(\Omega)}+\|g_{\phi}\|^{\alpha}_{L^{\alpha}(\Omega)})+
   C\int_{\Omega}d^{\alpha}(m_{0},\phi)dp
\end{align*}
for any $\mu$ where $m_{0}$ is a fixed point in $M$ and $g_{\phi}$
is a dominant function of $\phi$. The last inequalities show that
$\int_{\Omega}d^{\alpha}(u_{\mu}(p),m_{0})dp$ is uniformly
bounded. So we have
$$
\sup_{\mu}\{\int_{\Omega}d^{\alpha}(u_{\mu}(p),m_{0})dp+\|g_{\mu}\|^{\alpha}_{L^{\alpha}(\Omega)}\}\leq
C
$$
for a constant $C>0$ depending on $\phi$. Now we can use Theorem
\ref{thm:campactness} to get a subsequence $\{u_{\mu^{\prime}}\}$
of $\{u_{\mu}\}$ and $u\in{R^{1,\alpha}(\Omega,M)}$ such that
\begin{equation}\label{strong4}
   \lim_{\mu^{\prime}\rightarrow\infty}\int_{\Omega}d^{\alpha}(u_{\mu^{\prime}}(p),u(p))dp=0
\end{equation}
and
\begin{equation}\label{weakly4}
X_{l}u_{\mu}^{i}\textrm{ converges weakly in
}L^{\alpha}(\Omega)\textrm{ to }X_{l}u^{i}\textrm{ for }
i=1,\cdots,\widetilde{k}\textrm{ and }l=1,\cdots,k.
\end{equation}
From \eqref{strong4} and Theorem \ref{tracecharacterization} we
have $T(u)=T(\phi)$. Thus $u\in{R^{1}_{\phi}(\Omega,M)}$. From the
lower semicontinuity of $HE^{\alpha}(u,\Omega)$ with respect to
weak convergence and \eqref{weakly4}, we conclude that $u$ is a
minimizer.
\end{proof}

\section{Some remarks on the regularity of minimizers: Heisenberg group target}\label{sec:5}
In this section we briefly mention the known results to the
regularity problem. The regularity problem is still quite open and
new methods and tools should be developed to tackle it.\vskip10pt

For the case when the domain space is a C-C space and the target
is Euclidean, the H\"{o}lder regularities were obtain in
\cite{Ha1} and in \cite{JX} using different methods.

If the target
$M=(R^{\widetilde{n}},\widetilde{\Delta},\widetilde{d}_{c})$ is a
C-C space ($\widetilde{\Delta}$ is non-integrable), then since
$R^{1,\alpha}(\Omega,M)$ is not a linear space (because of the
contact condition), it is not trivial to construct contact
variations of the minimizer to deduce Euler-Lagrangian equations.
The simple example is the case studied by Capogna and Lin in
\cite{CL} where $\Omega\subset R^{n}$ is an Euclidean smooth
bounded domain and $M$ is the Heisenberg group $H^{m}$ with a
homogeneous metric $\rho$. We denote by $u=(z,t)=(x,y,t)$ elements
in $R^{1,\alpha}(\Omega,H^{m})$ ($\alpha\geq1$). Then from
Definition \ref{w1}, Theorem \ref{sec2:characterization} and
\eqref{item3} of Remark \ref{remark:cha} we have
$u\in{R^{1,\alpha}(\Omega,H^{m})}$ if and only if
\begin{enumerate}
    \item\label{H1} $z\in{W^{1,\alpha}(\Omega,R^{2m})}$;
    \item\label{H2} $t\in{L^{\frac{\alpha}{2}}}\cap \textrm{ACL}(\Omega,R)$;
    \item\label{H3} $u$ satisfies the Legendrian condition, that is, $\partial_{p_{i}}
t=2(y\cdot\partial_{p_{i}}x-x\cdot\partial_{p_{i}}y)$ a.e.
$p\in{\Omega}$ for $i=1,\cdots,n$.
\end{enumerate}
Here and in the sequel we denote by $x_{i}$ or $\partial_{p_{i}}x$
the partial derivative $\frac{\partial x}{\partial p_{i}}$ and
$\cdot$ denotes the inner product in $R^{n}$. Note that if $u$
satisfies \eqref{H1}, \eqref{H2} and \eqref{H3}, then from Sobolev
inequality and H\"{o}lder inequality,
$\partial_{p_{i}}t\in{L^{\beta}(\Omega)}$
($\beta=\frac{n\alpha}{2n-\alpha}$) automatically holds for
$i=1,\cdots,n$. Moreover if $\alpha\geq2$, then
$t\in{W^{1,\beta}}(\Omega)$. The horizontal energy is
$$
HE^{\alpha}(u,\Omega)=\int_{\Omega}|\nabla z(p)|^{\alpha}dp.
$$
\begin{lemma}\label{componenttrace}
Let $\alpha\geq2$ and
$u=(z_{u},t_{u}),v=(z_{v},t_{v})\in{R^{1,\alpha}(\Omega,H^{m})}$.
Then $T(u)=T(v)$ if and only if $T(z_{u})=T(z_{v})$ and
$T(t_u)=T(t_v)$ where $T(u)$ denotes the trace of $u$ on
$\partial\Omega$.
\end{lemma}
\begin{proof} The proof is straightforward. Since any homogeneous metrics
are equivalent,
\begin{align}
& C_{1}\left(|z_{u}-z_{v}|^{\alpha}+|t_{u}-t_{v}-2\omega(z_{u},z_{v})|^{\frac{\alpha}{2}}\right)\label{H4}\\
&\leq C_{2}\|u^{-1}v\|^{\alpha}\leq\rho^{\alpha}(u,v)\leq C_{3}\|u^{-1}v\|^{\alpha}\\
&\leq
C_{4}\left(|z_{u}-z_{v}|^{\alpha}+|t_{u}-t_{v}-2\omega(z_{u},z_{v})|^{\frac{\alpha}{2}}\right)\label{H5}
\end{align}
where $\|\cdot\|$ is the gauge norm:
$\|(z,t)\|=(|z|^{4}+t^{2})^{\frac{1}{4}}$; $\omega(\cdot,\cdot)$
is the standard symplectic form in $R^{2n}$, and $C_{i},i=1,2,3,4$
are constants. From Theorem \ref{tracecharacterization},
$T(u)=T(v)$ if and only if
$\rho(u,v)\in{W_{0}^{1,\alpha}(\Omega)}$. Since $\alpha\geq2$, the
term in \eqref{H4} belongs to $W^{1,1}(\Omega)$. Thus the
statement follows from \eqref{H4}-\eqref{H5} and the fact that
$\omega(z_{u},z_{v})|_{\partial\Omega}=0$ a.e if
$T(z_{u})=T(z_{v})$.
\end{proof}
The following lemma tells us that the projection of a weakly
contact map $u:\Omega\rightarrow H^{m}=R^{2m}\times R$ to $R^{2m}$
is a weakly isotropic map and conversely any weakly isotropic map
$z:\Omega\rightarrow R^{2m}$ can be lifted to be a weakly contact
map.
\begin{lemma}\label{iso-Le}
Let $\alpha\geq2$ and $\Omega\Subset U\Subset R^{n}$. If
$u=(z,t)=(x,y,t)\in{R^{1,\alpha}(\Omega,H^{m})}$, then
$z=(x,y)\in{W^{1,\alpha}(\Omega,R^{2m})}$ and satisfies the
following weakly isotropic condition:
$$
   z^{*}(\omega)=0 \textrm{ a.e in }\Omega,
$$
that is,
\begin{equation}\label{isotropic}
    x_{i}\cdot y_{j}=x_{j}\cdot y_{i} \textrm{ a.e. }p\in{\Omega}\textrm{ for
    }i,j=1,\cdots,n.
\end{equation}
 Conversely if
$\phi=(z_{\phi},t_{\phi})\in{R^{1,\alpha}(U,H^{m})}$ and
$$
z\in{W_{z_{\phi}}^{1,\alpha}(\Omega,R^{2m})}=:\{\widetilde{z}\in{W^{1,\alpha}(\Omega,R^{2m})}:T(\widetilde{z})=T(z_{\phi})\}
$$
and satisfies \eqref{isotropic},
 then there exists $t\in{W^{1,\beta}(\Omega)}$ such that
$u=(z,t)\in{R^{1,\alpha}(\Omega,H^{m})}$ and $t=t_{\phi}$ at
$\partial \Omega$.
\end{lemma}
\begin{proof} The first statement follows (essentially) from the
inequality (2.12) in \cite{Da}, see Lemma 2.12 and Theorem 2.16 in
\cite{CL} for details.

We prove the second statement. Let
$z=(x,y)\in{W_{z_{\phi}}^{1,\alpha}(\Omega,R^{2m})}$ and satisfies
\eqref{isotropic}. Let $\eta$ be the primitive form of the
standard symplectic form $\omega$ in $R^{2n}$, that is,
$d\eta=\omega$. We first prove that the 1-form
$$
\zeta(p)=z^{*}(\eta)=\frac{1}{2}\sum_{i=1}^{n}(x\cdot y_{i}-y\cdot
x_{i})dp_{i}
$$
belongs to $L^{\beta}(\Omega)$ and satisfies
\begin{equation}\label{exterior}
    d\zeta=0
\end{equation}
in the sense of distribution. It suffices to prove
\begin{equation}\label{exterior2}
    \int_{\Omega}(x\cdot y_{i}-y\cdot x_{i})\varphi_{j}-(x\cdot y_{j}-y\cdot
    x_{j})\varphi_{i}dp=0
\end{equation}
for any $\varphi\in{C_{0}^{\infty}}$ and $i\neq j$. We mollify $z$
and let
$z^{\epsilon}=(x^{\epsilon},y^{\epsilon})=z*\delta_{\epsilon}$
where $\delta_{\epsilon}$ is a standard mollifier. We have
\begin{align*}
   \textrm{I}&= \int_{\Omega}(x^{\epsilon}\cdot y^{\epsilon}_{i}-y^{\epsilon}\cdot x^{\epsilon}_{i})\varphi_{j}-
    (x^{\epsilon}\cdot y^{\epsilon}_{j}-y^{\epsilon}\cdot
    x^{\epsilon}_{j})\varphi_{i}dp\\
    &=- \int_{\Omega}(x^{\epsilon}\cdot y^{\epsilon}_{i}-y^{\epsilon}\cdot x^{\epsilon}_{i})_{j}\varphi-
    (x^{\epsilon}\cdot y^{\epsilon}_{j}-y^{\epsilon}\cdot
    x^{\epsilon}_{j})_{i}\varphi dp\\
    &=\textrm{II}+\textrm{III}
\end{align*}
where
$$
\textrm{II}=\int_{\Omega}(x_{i}^{\epsilon}\cdot
y_{j}^{\epsilon}-y^{\epsilon}_{i}\cdot
x_{j}^{\epsilon})\varphi-(x_{j}^{\epsilon}\cdot
y_{i}^{\epsilon}-y^{\epsilon}_{j}\cdot x_{i}^{\epsilon})\varphi dp
$$
and
\begin{align*}
\textrm{III}&=\int_{\Omega}x^{\epsilon}\cdot
\left((y^{\epsilon}_{j})_{i}-(y^{\epsilon}_{i})_{j}\right)\varphi+y^{\epsilon}\cdot
\left((x_{i})^{\epsilon}_{j}-(x_{j})^{\epsilon}_{i}\right)\varphi
dp\\
&=0.
\end{align*}
 Since $z$ satisfies \eqref{isotropic}, we have
$ \lim_{\epsilon\rightarrow0}\textrm{II}=0. $ Consequently
\eqref{exterior2} follows and \eqref{exterior} holds. Thus we can
apply
 a well known result about boundary
value problem involving differential forms, see e.g. Chapter 3 of
\cite{Sc}, to get a $t\in{W^{1,\beta}(\Omega)}$ of the following
equation
\begin{equation}\label{equation}
\left\{%
\begin{array}{ll}
   dt=4\zeta, & \hbox{in $\Omega$;} \\
    t=t_{\phi}, & \hbox{in $\partial\Omega$.} \\
\end{array}%
\right.
\end{equation}
Note that \eqref{equation} means $\partial_{p_{i}}
t=2(y\partial_{p_{i}}x-x\partial_{p_{i}}y)$ a.e. $p\in{\Omega}$
for $i=1,\cdots,n$.
\end{proof}

Let
$$
W_{I}^{1,\alpha}(\Omega,R^{2m})=:\{z\in{W^{1,\alpha}(\Omega,R^{2m})}:z
\textrm{ satisfies }\eqref{isotropic}\}
$$
 be the Sobolev space of weakly isotropic mappings. If
 $\varphi\in{W_{I}^{1,\alpha}(\Omega,R^{2m})}$,
 we denote by
$$
W_{I,\varphi}^{1,\alpha}(\Omega,R^{2m})=:W_{I}^{1,\alpha}(\Omega,R^{2m})\cap
W_{\varphi}^{1,\alpha}(\Omega,R^{2m})
$$
the Sobolev space of weakly isotropic mappings with the trace of
$\varphi$.

From Lemma \ref{componenttrace} and \ref{iso-Le}, we get the
following theorem.
\begin{theorem}\label{deduction}
Let $\alpha\geq2$ and $\Omega\Subset U\Subset R^{n}$ be two
bounded open sets. Let
$\phi=(z_{\phi},t_{\phi})\in{R^{1,\alpha}(U,H^{m})}$. Then
$u=(z,t)\in{R_{\phi}^{1,\alpha}(\Omega,H^{m})}$ is a solution of
Problem \eqref{dirichletproblem}, that is,
$$
HE^{\alpha}(u,\Omega)=\inf_{v\in{R_{\phi}^{1,\alpha}(\Omega,M)}}HE^{\alpha}(v,\Omega),
$$
 if and only if $z$
is a solution of the following isotropically constrained
variational problem:
\begin{equation}\label{isotropicconstrained}
\textrm{ to find
}z_{0}\in{W_{I,z_{\phi}}^{1,\alpha}(\Omega,R^{2m})} \textrm{ such
that
}E(z_{0},\Omega)=\inf_{z\in{W_{I,z_{\phi}}^{1,\alpha}(\Omega,R^{2m})}}E(z,\Omega)
\end{equation}
where $E(z,\Omega)=\int_{\Omega}|\nabla z(p)|^{\alpha}dp$.
\end{theorem}
\begin{remark}\((1)\) The existence of solutions to Problem
\eqref{isotropicconstrained} can be easily established.

\((2)\) When $n=\alpha=2$, due to the conformal invariance of the
Dirichlet integral, Problem \eqref{isotropicconstrained} is
closely related to the following isotropically constrained Plateau
problem studied by \cite{SW} and \cite{Qiu}:
\begin{equation}\label{plateau}
    \textrm{ to find }l_{0}\in{\mathcal{X}_{I,\Gamma}}\textrm{ such that
    }Area(l_{0},B_{1})=\inf_{l\in{\mathcal{X}_{I,\Gamma}}}Area(l,B_{1})
\end{equation}
where $\Gamma$ is a piecewise $C^{1}$ closed Jordan curve such
that $\int_{\Gamma}\eta=0$; $B_{1}$ is the ball centered at 0 with
radius 1;
$\mathcal{X}_{I,\Gamma}=\{l\in{W^{1,2}_{I}(B_{1},R^{2m})}:l_{\partial
B_{1}}$ is continuous and is a monotone map onto $\Gamma$\};
$Area(l,B_{1})$ is the area of the image of $l$.
\end{remark}
\begin{proposition}\label{lipschitz}
Let $\Omega\subset R^{2}$ be a bounded open set and let
$\alpha\geq 2$. Let $\phi\in{R^{1,\alpha}(\Omega,H^{1})}$. If
$u\in{R^{1,\alpha}(\Omega,H^{1})}$ such that
$$
HE^{\alpha}(u,\Omega)=\inf_{v\in{R_{\phi}^{1,\alpha}(\Omega,H^{1})}}HE^{\alpha}(v,\Omega),
$$
then $u$ is Lipschitz continuous (with respect to the C-C metric)
in the interior of $\Omega$.
\end{proposition}
Proposition \ref{lipschitz} was proven in Theorem 4.3 of \cite{CL}
where the authors asserted that Proposition \ref{lipschitz} should
hold for $m\geq1$. In our opinion, since the equation (4.2) in
\cite{CL} holds only for $m=1$, the case when $m>1$ remains
open.\vskip10pt

Keeping Theorem \ref{deduction} in mind we see the following
theorem is just a copy of the main result in \cite{Qiu}.
\begin{theorem}
Let $U$ be a neighborhood of  a smooth, connected and simply
connected, bounded open set $\Omega$ such that $\Omega\Subset
U\Subset R^{2}$ and let $m>1$. Let
$\phi=(z_{\phi},t_{\phi})\in{R^{1,2}(U,H^{m})}$ such that
$\phi|_{\partial\Omega}$ is continuous and monotone onto
$\phi(\partial\Omega)$ which is a piecewise $C^{1}$ closed Jordan
curve in $R^{2m}$. Let $u=(z,t)\in{R^{1,2}_{\phi}(\Omega,H^{m})}$
be a solution of Problem \eqref{dirichletproblem}, that is
$$
HE^{2}(u,\Omega)=\inf_{v\in{R_{\phi}^{1,2}(\Omega,M)}}HE^{2}(v,\Omega).
$$
Then $z$ is H\"{o}lder continuous in $\Omega$ and $u$ is smooth in
$\Omega$ with possibly isolated singularities.
\end{theorem}

\end{document}